# Robust Transceiver Optimization for Downlink Coordinated Base Station Systems: Distributed Algorithm .

Tadilo Endeshaw Bogale, *Student Member, IEEE*, Luc Vandendorpe, *Fellow, IEEE* and
Batu Krishna Chalise, *Member, IEEE*
ICTEAM Institute
Universitè catholique de Louvain
Place du Levant, 2, B-1348, Louvain La Neuve, Belgium
Email: {tadilo.bogale, luc.vandendorpe}@uclouvain.be and batu.chalise@villanova.edu

*Abstract*— This paper considers the joint transceiver design for downlink multiuser multiple-input single-output (MISO) systems with coordinated base stations (BSs) where imperfect channel state information (CSI) is available at the BSs and mobile stations (MSs). By incorporating antenna correlation at the BSs and taking channel estimation errors into account, we solve two robust design problems: 1) minimizing the weighted sum of mean-square-error (MSE) with per BS antenna power constraint, and 2) minimizing the total power of all BSs with per user MSE target and per BS antenna power constraints. These problems are solved as follows. First, for fixed receivers, we propose centralized and novel computationally efficient distributed algorithms to jointly optimize the precoders of all users. Our centralized algorithms employ the second-order-cone programming (SOCP) approach, whereas, our novel distributed algorithms use the Lagrangian dual decomposition, modified matrix fractional minimization and an iterative method. Second, for fixed BS precoders, the receivers are updated by the minimum mean-square-error (MMSE) criterion. These two steps are repeated until convergence is achieved. In all of our simulation results, we have observed that the proposed distributed algorithms achieve the same performance as that of the centralized algorithms. Moreover, computer simulations verify the robustness of the proposed robust designs compared to the non-robust/naive designs.

## I. INTRODUCTION

The next generation multimedia communications are expected to support high data rates. To meet this demand multi-antenna systems are recommended as they significantly increase the spectral efficiency of wireless channels [1], [2]. The performance enhancement is achieved by exploiting the transmit and receive diversity. In [3], a fundamental relation between mutual information and minimum mean-square-error (MMSE) has been established for multiple-input multiple-output (MIMO) Gaussian channels. Furthermore, it has been shown that different transceiver optimization problems are equivalently reformulated as a function of MMSE matrix, for instance, minimizing bit error rate, maximizing capacity etc [4]–[6]. For these reasons, mean-square-error (MSE)-based design problems are commonly examined in multiuser networks.

In general, the uplink channel MSE-based problems are better understood than that of the downlink channel problems. For this reason, most of the research works examine MSE-based problems for the downlink multiuser systems [5], [7]–[12]. However, all of the these papers examine their problems for conventional downlink networks. In these networks, base stations (BSs) from different cells communicate with their respective remote terminals independently. Hence, in the latter network, inter-cell interference is obliged to be considered as a background noise. Recently, it has been shown that BS coordination communication is a promising technique to significantly improve the capacity of wireless channels by mitigating (or possibly canceling) inter-cell interference [13]–[15]. The BS coordination can be performed by two approaches. In the first approach, BSs are coordinated at the beamforming (precoder) level [14], whereas in the second approach, coordination takes place both at the signal and beamforming (precoder) levels [13], [15]. It is well know that the latter coordination approach has better performance gain compared to the former one [15], [16]. This performance improvement, however, requires additional signal coordination. In the current paper, we focus on the second BS coordination approach (the approach of [13] and [15]). In [17], four MSE-based linear transceiver optimization problems have been considered for coordinated BS MIMO systems. These problems are examined by assuming that the total power of each BS or the individual power of each BS antenna is constrained. The optimization problems of [17] are solved as follows. First, by keeping the receivers constant, the precoders of all users are jointly optimized using a second-order-cone programming (SOCP) approach (SOCP problems are convex and can be solved using interior point (IP) methods [18]). Second, for the given

The authors would like to thank the Region Wallonne for the financial support of the project MIMOCOM in the framework of which this work has been achieved. Part of this work has been published in the 44th Annual Asilomar Conference on Signals, Systems, and Computers, Pacific Grove, California, USA, Nov. 2010

Tadilo E. Bogale and Luc Vandendorpe are with the ICTEAM Institute, Université catholique de Louvain, Place du Levant 2, 1348 - Louvain La Neuve, Belgium. Email: {tadilo.bogale, luc.vandendorpe}@uclouvain.be, Phone: +3210478071, Fax: +3210472089. Batu K. Chalise was with the ICTEAM Institute, Universite catholique de Louvain, Belgium when this work was completed. He is currently with the Center for Advanced Communications, Villanova University, 800 Lancaster Avenue, Villanova, PA, 19085, USA, Phone: +16105197371, Fax: +16105196118, Email: batu.chalise@villanova.edu.

BS precoders, the receiver of each user is optimized by the MMSE method. The first and second steps are repeated in an iterative manner to jointly optimize the transmitters and receivers. Thus, in [17], the receiver of each user can be optimized independently and distributively. However, the joint optimization of the precoders has been carried out by a centralized algorithm. When the number of users and/or BSs increase, the computational cost of the joint precoder design also increases [19]. Consequently, solving the precoder optimization problem in a centralized manner, especially for large-scale coordinated networks, is not a computationally efficient approach. This motivates us to develop distributed algorithms to solve MSE-based problems for coordinated BS systems with per BS antenna power constraints in [20]. This paper solves its optimization problems distributively by applying the Lagrangian dual decomposition, modified matrix fractional minimization and an iterative technique.

In the current paper, we extend the work of [20] to robust case. The goal of this work is to jointly optimize the transmitters and receivers of all users when imperfect channel state information (CSI) is available both at the BSs and mobile stations (MSs), and with antenna correlation at the BSs. It is known in [21]–[23] that transmit antenna correlation matrices depend on array parameters (such as array geometry, antenna spacing), and the average angle of arrival (AOA) of scattered signals from user and the corresponding angular spread. This means that the transmit antenna correlation matrices (which capture spatial variation) vary at a rate much slower than the fast fading component (that captures temporal variations) of downlink channel. Thus, errors caused from the estimation of fast fading part of the channel can significantly outnumber the errors caused from the estimation of slowly varying antenna correlation matrices. Based upon these discussions, it is clear that the transmit correlation matrices can be obtained from long term channel statistics with a reasonable accuracy [24]. Due to this reason, perfect transmit antenna correlation matrices are assumed to be available at the BSs and MSs. We also assume that each BS is equipped with multiple antennas and each MS is equipped with single antenna. For our robust transceiver designs, a stochastic approach has been utilized. We design the transmitters and receivers by considering the practically relevant scenario where spatial correlation matrices as seen from each BS are different for different MSs and the variance of the estimation errors corresponding to the estimates of the channels are different. For this CSI model, we examine the following MSE-based robust design problems.

1) The robust minimization of the weighted sum MSE with per BS antenna power constraint ($\mathcal{P}1$).
2) The robust minimization of the total sum power of all BSs with per BS antenna power and per user MSE target constraints ($\mathcal{P}2$)[1].

To the best of our knowledge, problems $\mathcal{P}1$ and $\mathcal{P}2$ are non-convex. Hence, convex optimization tools can not be used to solve them. Each of these problems are solved iteratively as follows. First, for given receivers, we propose centralized and novel computationally efficient distributed algorithms to design the optimal precoders of all users. Our centralized algorithm designs the precoders of all users using SOCP approach and our novel distributed algorithm designs the precoders of all users by employing the Lagrangian dual decomposition, modified matrix fractional minimization and an iterative approach. Second, like in [17] and [20], the receiver of each user is optimized independently using minimum average mean-square-error (MAMSE) approach. These steps are repeated until convergence is achieved. The centralized and distributed algorithms require the complete channel estimates of all MSs. The centralized precoder design algorithms are developed by extending the approach of [17] to the case where imperfect CSI is available at the BSs and MSs. However, this extension does not change the fact that the robust problem can be reformulated into a SOCP problem as in the case of non-robust problem [17]. Thus, the novelty of our current paper relies mainly on the proposed distributed algorithms where we have used modified matrix fractional minimization techniques to solve the precoder design problems of $\mathcal{P}1$ and $\mathcal{P}2$ distributively. As a result, our new distributed algorithms are able to solve the transceiver design problems of $\mathcal{P}1$ and $\mathcal{P}2$ with less computational cost than that of the centralized algorithms. Furthermore, in all of our simulation results, we have observed that our novel distributed algorithms achieve the same performance as that of the centralized algorithms. The main contributions of this paper is thus summarized as follows.

1) We propose novel computationally efficient distributed algorithms to jointly optimize the precoders of all users for problems $\mathcal{P}1$ and $\mathcal{P}2$. As will be clear later, the proposed distributed algorithms can be extended straightforwardly to solve $\mathcal{P}1$ and $\mathcal{P}2$ for MIMO coordinated BS systems with per BS antenna (groups of BS antennas) power constraints.
2) We have demonstrated that the proposed distributed algorithms for $\mathcal{P}1$ and $\mathcal{P}2$ achieve the same performance as that of the centralized algorithms proposed for $\mathcal{P}1$ and $\mathcal{P}2$.
3) We examine the joint effect of channel estimation errors and antenna correlations on the performance of $\mathcal{P}1$ and $\mathcal{P}2$.

The remaining part of this paper is organized as follows. We present the coordinated multi-antenna BS system model in Section II. In Section III, the multiuser channel model under imperfect CSI is presented. Section IV discusses the robust design problems $\mathcal{P}1$ - $\mathcal{P}2$, and the proposed centralized and distributed algorithms. The extensions of our centralized and distributed algorithms for $\mathcal{P}1$ and $\mathcal{P}2$ in a MIMO coordinated BS system is discussed in Section V. In Section VI, computer simulations are used to compare the performance of the centralized and distributed algorithms, and the robust and non-robust/naive designs. Finally, conclusions are drawn in Section VII.

---

[1]In any transceiver design problem, the notion of power arises at the transmitter side. In this paper, since we are examining downlink transceiver design problems, the power constraint appears only at the BSs. In a practical downlink multi-antenna BS systems, each BS antenna has its own power amplifier and the maximum power of each antenna is limited [25]. This motivates us to consider the power constraint of each BS antenna. However, as will be clear later, the proposed algorithms for $\mathcal{P}1$ and $\mathcal{P}2$ can be extended straightforwardly to handle the sum power constraint of the whole network or groups of antennas.

*Notations:* The following notations are used throughout this paper. Upper/lower case boldface letters denote matrices/column vectors. $\text{vec}(\mathbf{X})$, $\text{tr}(\mathbf{X})$, $\mathbf{X}^\star$, $\mathbf{X}^T$, $\mathbf{X}^H$ and $\text{E}(\mathbf{X})$ denote vectorization, trace, optimal, transpose, conjugate transpose and expected value of $\mathbf{X}$, respectively. $\mathbf{I}_n(\mathbf{I})$ is an identity matrix of size $n \times n$ (appropriate size) and $\mathcal{C}^{M\times M}$ represent spaces of $M \times M$ matrices with complex entries. The diagonal and block-diagonal matrices are represented by $\text{diag}(.)$ and $\text{blkdiag}(.)$ respectively. Subject to is denoted by s.t, $\|\mathbf{x}\|_n$ is the the $n$th norm of a vector $\mathbf{x}$.

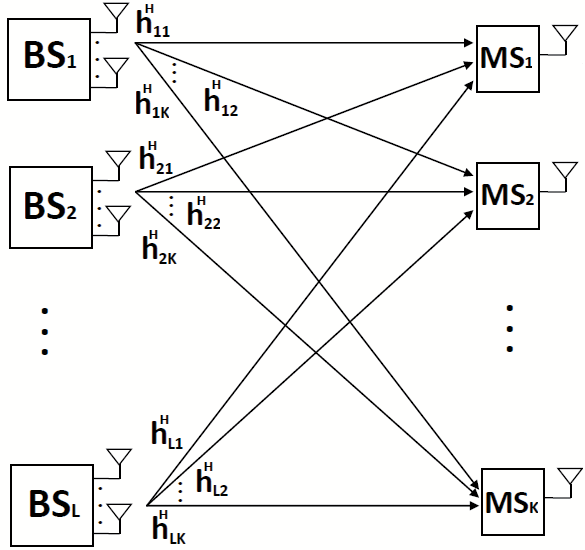

Fig. 1. Coordinated base station system model.

## II. SYSTEM MODEL

We consider a coordinated BS system as shown in Fig. 1 where $L$ BSs are serving $K$ decentralized single antenna MSs. The $l$th BS is equipped with $N_l$ transmit antennas. By denoting the symbol intended for the $k$th user as $d_k$, the entire symbol can be written in a data vector $\mathbf{d} \in \mathcal{C}^{K\times 1}$ as $\mathbf{d} = [d_1,\cdots,d_K]^T$. The $l$th BS precodes $\mathbf{d}$ into an $N_l$ length vector by using its overall precoder matrix $\mathbf{B}_l = [\mathbf{b}_{l1},\cdots,\mathbf{b}_{lK}]$, where $\mathbf{b}_{lk} \in C^{N_l \times 1}$ is the precoder vector of the $l$th BS for the $k$th MS. For convenience, we follow the same channel vector notations as in [20]. The $k$th MS employs a receiver $w_k \in C$ to estimate its symbol $d_k$. The estimated symbol at the $k$th MS is given by

$$\hat{d}_k = w_k^H(\sum_{l=1}^{L}\mathbf{h}_{lk}^H\mathbf{B}_l\mathbf{d} + n_k) = w_k^H(\mathbf{h}_k^H\sum_{i=1}^{K}\mathbf{b}_i d_i + n_k) \quad (1)$$

where $\mathbf{h}_k^H = [\mathbf{h}_{1k}^H,\cdots,\mathbf{h}_{Lk}^H] \in \mathcal{C}^{1\times N}$, $\mathbf{B} = [\mathbf{B}_1;\cdots;\mathbf{B}_L]$, $\mathbf{b}_k = [\mathbf{b}_{1k}^T,\cdots,\mathbf{b}_{Lk}^T]^T \in \mathcal{C}^{N\times 1}$ is the precoder vector for the $k$th user, $N = \sum_{l=1}^{L}N_l$, $\mathbf{h}_{lk}^H \in \mathcal{C}^{1\times N_l}$ is the channel vector between the $l$th BS and the $k$th MS, and $n_k$ is the additive noise at the $k$th MS. It is clearly seen that the last expression of (1) has exactly the same form as the estimate of $d_k$ for the downlink MISO system where a BS equipped with $N$ transmitt antennas is serving $K$ decentralized users. Hence, we can interpret coordinated BS system as a one giant downlink system [17], [19]. It is assumed that $n_k$ is a zero-mean circularly symmetric complex Gaussian (ZMCSCG) random variable with the variance $\sigma_k^2$, i.e., $n_k \sim \mathcal{NC}(0,\sigma_k^2)$. We also assume that the symbol $d_k$ is a ZMCSCG random variable with unit variance and is independent of $\{d_i\}_{i=1,i\neq k}^{K}$ and noise $n_k$, i.e., $\text{E}\{d_k d_k^H\} = 1$, $\text{E}\{d_k d_i^H\} = 0, \forall i \neq k$ and $\text{E}\{d_k n_k^H\} = 0$. For this system model, when perfect CSI is available at the BSs and MSs, the MSE of the $k$th user can be expressed as

$$\xi_k = \text{E}_{\mathbf{d},n_k}\{(\hat{d}_k - d_k)(\hat{d}_k - d_k)^H\}$$
$$= w_k^H\left(\left[\sum_{l=1}^{L}\mathbf{h}_{lk}^H\mathbf{B}_l\right]\left[\sum_{l=1}^{L}\mathbf{h}_{lk}^H\mathbf{B}_l\right]^H + \sigma_k^2\right)w_k -$$
$$w_k^H\sum_{l=1}^{L}\mathbf{h}_{lk}^H\mathbf{b}_{lk} - \sum_{l=1}^{L}\mathbf{b}_{lk}^H\mathbf{h}_{lk}w_k + 1.$$

## III. CHANNEL MODEL

Considering antenna correlation at the BSs, we model the Rayleigh fading MISO channel between the $l$th BS and the $k$th MS as $\mathbf{h}_{lk}^H = \mathbf{h}_{wlk}^H\tilde{\mathbf{R}}_{blk}^{1/2}$, where the elements of $\mathbf{h}_{wlk}^H$ are independent and identically distributed (i.i.d) ZMCSCG random variables all with unit variance and $\tilde{\mathbf{R}}_{blk} \in C^{N_l\times N_l}$ is the antenna correlation matrix as seen from the $l$th BS [26], [24]. The channel estimation of the $k$th MS ($\mathbf{h}_{lk}^H$) can be performed on $\mathbf{h}_{wlk}^H, \forall l$, using an orthogonal training method [27]. Upon doing so, the true channel between the $l$th BS and $k$th MS $\mathbf{h}_{lk}^H$ is given by [12]

$$\mathbf{h}_{lk}^H = (\widehat{\mathbf{h}}_{wlk}^H + \mathbf{e}_{wlk}^H)\tilde{\mathbf{R}}_{blk}^{1/2} = \widehat{\mathbf{h}}_{lk}^H + \mathbf{e}_{wlk}^H\tilde{\mathbf{R}}_{blk}^{1/2} \quad (2)$$

where $\widehat{\mathbf{h}}_{wlk}^H$ is the MMSE estimate of $\mathbf{h}_{wlk}^H$, $\widehat{\mathbf{h}}_{lk}^H = \widehat{\mathbf{h}}_{wlk}^H\tilde{\mathbf{R}}_{blk}^{1/2}$ and $\mathbf{e}_{wlk}^H$ is the estimation error for which its entries are i.i.d with $C\mathcal{N}(0,\sigma_{elk}^2)$. In the case of linearly dependent channel estimation errors, $\mathbf{e}_{wlk}^H$ and $\mathbf{h}_{lk}^H$ can be expressed as $\mathbf{e}_{wlk}^H = \tilde{\mathbf{e}}_{wlk}^H\mathbf{Z}_{lk}$ and $\mathbf{h}_{lk}^H = \widehat{\mathbf{h}}_{lk}^H + \tilde{\mathbf{e}}_{wlk}^H\bar{\tilde{\mathbf{R}}}_{blk}$, where $\mathbf{Z}_{lk} \in C^{\tilde{N}_{lk}\times N_l}$, $\tilde{N}_{lk} \leq N_l$, $\tilde{\mathbf{e}}_{wlk}^H \in C^{1\times \tilde{N}_{lk}} = C\mathcal{N}(0,\tilde{\sigma}_{elk}^2)$ and $\bar{\tilde{\mathbf{R}}}_{blk} = \mathbf{Z}_{lk}\tilde{\mathbf{R}}_{blk}^{1/2}$. For simplicity, the current paper examines $\mathcal{P}1$ and $\mathcal{P}2$ for the channel model given in (2). As will be clear later, the approaches of this paper can also be applied to solve $\mathcal{P}1$ and $\mathcal{P}2$ with linearly dependent channel estimation error models.

The main idea of the robust design is that $\{\mathbf{e}_{wlk}^H\}_{k=1}^{K}, \forall l$ are unknown but $\{\widehat{\mathbf{h}}_{wlk}^H, \tilde{\mathbf{R}}_{blk}$ and $\sigma_{elk}^2\}_{k=1}^{K}, \forall l$ are available. We assume that the $k$th MS estimates its channel (i.e., $\widehat{\mathbf{h}}_{lk}^H, \forall l$) and feeds $\{\widehat{\mathbf{h}}_{lk}^H \text{ and } \sigma_{elk}^2\}, \forall l$ back to the BSs without any error and delay. Since the $l$th BS has the channel estimates of all MSs, it can compute $\widetilde{\mathbf{R}}_{blk}$ locally from the long term channel statistics of $\widehat{\mathbf{h}}_{lk}^H$ [24]. Thus, both the BSs and MSs have the same channel imperfections. The average mean-square-error (AMSE) of the $k$th MS ($\bar{\xi}_k$) can be expressed as

$$\begin{aligned}
\bar{\xi}_k =& \mathrm{E}_{\mathbf{e}_{wk}^H}\{\xi_k\} \\
=& w_k^H \bigg(\bigg[\sum_{l=1}^L \widehat{\mathbf{h}}_{lk}^H \mathbf{B}_l\bigg]\bigg[\sum_{l=1}^L \widehat{\mathbf{h}}_{lk}^H \mathbf{B}_l\bigg]^H + \\
& \sum_{l=1}^L \sigma_{elk}^2 \mathrm{tr}\{\widetilde{\mathbf{R}}_{blk}^{1/2} \mathbf{B}_l \mathbf{B}_l^H \widetilde{\mathbf{R}}_{blk}^{1/2}\} + \sigma_k^2\bigg)w_k - \\
& w_k^H \sum_{l=1}^L \widehat{\mathbf{h}}_{lk}^H \mathbf{b}_{lk} - \sum_{l=1}^L \mathbf{b}_{lk}^H \widehat{\mathbf{h}}_{lk} w_k + 1 \\
=& w_k^H (\widehat{\mathbf{h}}_k^H \sum_{i=1}^K \mathbf{b}_i \mathbf{b}_i^H \widehat{\mathbf{h}}_k + \mathrm{tr}\{\mathbf{R}_{bk}^{1/2} \sum_{i=1}^K \mathbf{b}_i \mathbf{b}_i^H \mathbf{R}_{bk}^{1/2}\} \\
& + \sigma_k^2) w_k - w_k^H \widehat{\mathbf{h}}_k^H \mathbf{b}_k - \mathbf{b}_k^H \widehat{\mathbf{h}}_k w_k + 1 \\
=& w_k^H (\widehat{\mathbf{h}}_k^H \sum_{i=1}^K \mathbf{b}_i \mathbf{b}_i^H \widehat{\mathbf{h}}_k + \sum_{i=1}^K \mathbf{b}_i^H \mathbf{R}_{bk} \mathbf{b}_i + \sigma_k^2) w_k \\
& - w_k^H \widehat{\mathbf{h}}_k^H \mathbf{b}_k - \mathbf{b}_k^H \widehat{\mathbf{h}}_k w_k + 1
\end{aligned} \quad (3)$$

where $\widehat{\mathbf{h}}_k^H = [\widehat{\mathbf{h}}_{1k}^H, \cdots, \widehat{\mathbf{h}}_{Lk}^H]$, $\mathbf{R}_{bk} = \mathrm{blkdiag}(\sigma_{e1k}^2 \widetilde{\mathbf{R}}_{b1k}, \cdots, \sigma_{eLk}^2 \widetilde{\mathbf{R}}_{bLk})$, and in the second equality we use the fact $\mathrm{E}_{\mathbf{e}}\{\mathbf{e}^H \boldsymbol{\Phi} \mathbf{e}\} = \sigma_e^2 \mathrm{tr}\{\boldsymbol{\Phi}\}$, if the entries of $\mathbf{e}$ are i.i.d with $\mathcal{CN}(0, \sigma_e^2)$ and $\boldsymbol{\Phi}$ is a given matrix [12], [28]. Note that one can extend the precoder/decoder design problems of $\mathcal{P}1$ and $\mathcal{P}2$ by incorporating channel estimation and feedback errors with time delay. In this case, the expression of $\{\bar{\xi}_k\}_{k=1}^K$ will be different from (3). Hence, solving $\mathcal{P}1$ and $\mathcal{P}2$ with erroneous feedback and non-zero delay is an open research topic.

## IV. Problem Formulations and proposed Solutions

In this section, we examine problems $\mathcal{P}1$ and $\mathcal{P}2$. For each problem, we use the following general optimization framework. First, for fixed receivers, the precoders of all users are optimized. Second, using the latter precoders, the optimal receiver of each user is designed by the MAMSE receiver method. Since designing the receivers by MAMSE approach is optimal for any robust MSE-based problem, $\mathcal{P}1$ and $\mathcal{P}2$ utilize the same MAMSE expression to design their receivers. Finally, these two steps are repeated until convergence is achieved. The MAMSE receiver of the $k$th user is given by [20]

$$w_k = \frac{\widehat{\mathbf{h}}_k^H \mathbf{b}_k}{\widehat{\mathbf{h}}_k^H \sum_{i=1}^K \mathbf{b}_i \mathbf{b}_i^H \widehat{\mathbf{h}}_k + \mathrm{tr}\{\mathbf{R}_{bk} \sum_{i=1}^K \mathbf{b}_i \mathbf{b}_i^H\} + \sigma_k^2}. \quad (4)$$

In the following, for fixed receivers $\{w_k\}_{k=1}^K$, we propose centralized and computationally efficient distributed precoder design algorithms for problems $\mathcal{P}1$ and $\mathcal{P}2$.

### A. Robust weighted sum MSE minimization problem ($\mathcal{P}1$)

The robust weighted sum MSE minimization with per BS antenna power constraint problem can be formulated as

$$\min_{\{w_k, \mathbf{b}_k\}_{k=1}^K} \sum_{k=1}^K \eta_k \bar{\xi}_k, \quad \text{s.t} \left[\sum_{k=1}^K \mathbf{b}_k \mathbf{b}_k^H\right]_{n,n} \leq p_n, \forall n \quad (5)$$

where $\eta_k$ is the AMSE weighting factor of the $k$th user and $p_n$ is the maximum available power at the $n$th BS antenna. The antenna numbers are assigned from the first antenna of $\mathrm{BS}_1$ (which corresponds to antenna 1) to the last antenna of $\mathrm{BS}_L$ (which corresponds to antenna $N$).

*1) Centralized precoder design for $\mathcal{P}1$:* The objective function of the above problem can be expressed as

$$\begin{aligned}
\sum_{k=1}^K \eta_k \bar{\xi}_k =& \sum_{k=1}^K \bigg\{ \mathbf{b}_k^H (\sum_{i=1}^K \eta_i \widehat{\mathbf{h}}_i w_i w_i^H \widehat{\mathbf{h}}_i^H + \eta_i w_i w_i^H \mathbf{R}_{bi}) \mathbf{b}_k \\
& - \eta_k w_k^H \widehat{\mathbf{h}}_k^H \mathbf{b}_k - \eta_k \mathbf{b}_k^H \widehat{\mathbf{h}}_k w_k + \sigma_k^2 \eta_k w_k^H w_k + \eta_k \bigg\} \\
=& \mathrm{tr}\{(\sqrt{\boldsymbol{\eta}} \mathbf{W}^H \widehat{\mathbf{H}}^H \mathbf{B} - \sqrt{\boldsymbol{\eta}})^H (\sqrt{\boldsymbol{\eta}} \mathbf{W}^H \widehat{\mathbf{H}}^H \mathbf{B} - \sqrt{\boldsymbol{\eta}})\} \\
& + \mathrm{tr}(\mathbf{B}^H \boldsymbol{\Psi} \mathbf{B}) + \mathrm{tr}(\boldsymbol{\sigma}^2 \boldsymbol{\eta} \mathbf{W}^H \mathbf{W})
\end{aligned} \quad (6)$$

where $\boldsymbol{\eta} = \mathrm{diag}(\eta_1, \cdots, \eta_K)$, $\mathbf{W} = \mathrm{diag}(w_1, \cdots, w_K)$, $\widehat{\mathbf{H}} = [\widehat{\mathbf{h}}_1, \cdots, \widehat{\mathbf{h}}_K]$, $\boldsymbol{\sigma} = \mathrm{diag}(\sigma_1, \cdots, \sigma_K)$ and $\boldsymbol{\Psi} = \sum_{i=1}^K \eta_i w_i w_i^H \mathbf{R}_{bi}$. For fixed receivers $\{w_k\}_{k=1}^K$ and using (6), the global optimal $\{\mathbf{b}_k\}_{k=1}^K$ of (5) can be obtained by solving the following problem [17]

$$\min_{\chi, \mathbf{B}} \chi \quad \text{s.t} \quad \|\boldsymbol{\mu}\|_2 \leq \chi, \quad \|\widetilde{\mathbf{b}}_n\|_2 \leq \sqrt{p_n}, \forall n \quad (7)$$

where $\boldsymbol{\mu} = [\mathrm{vec}(\sqrt{\boldsymbol{\eta}} \mathbf{W}^H \widehat{\mathbf{H}}^H \mathbf{B} - \sqrt{\boldsymbol{\eta}}) \, ; \, \mathrm{vec}(\sqrt{\boldsymbol{\Psi}} \mathbf{B})]$ and $\widetilde{\mathbf{b}}_n^H$ as the $n$th row of $\mathbf{B}$. As we can see, (7) is a SOCP problem for which the global optimal solution is obtained by existing convex optimization tools [18]. This problem has $2NK+1$ real optimization variables, $N$ second-order-cone (SOC) constraints where each of them consists of $2K$ real dimensions and one SOC constraint with $2K(K+N)$ real dimensions. According to [29] (see page 196), the computational complexity of the latter problem in terms of number of iterations is upper bounded by $O(\sqrt{N+1})$ where the complexity of each iteration is within the order of $O((2KN+1)^2(2K^2+4KN))$. Thus, the total worst-case computational complexity of (7) is given by $O(\sqrt{N+1}(2KN+1)^2(2K^2+4KN))$. This shows that for large networks, the centralized precoder design scheme appears to be impractical. This motivates us to design the precoders of each user distributively with less computational cost than that of the centralized precoder design approach.

*2) Distributed precoder design for $\mathcal{P}1$:* For fixed $\{w_k\}_{k=1}^K$, to solve the precoders of (5) distributively, we propose the Lagrangian dual decomposition technique[2]. To this end, we first express the Lagrangian function associated with

---
[2]Since the precoder design problem of (5) is convex, there is a zero duality gap between the primal and dual problems.

(5) as

$$L(\boldsymbol{\lambda}, \mathbf{B}) = \sum_{k=1}^{K} \eta_k \bar{\xi}_k + \sum_{n=1}^{N} \lambda_n \left( [\sum_{i=1}^{K} \mathbf{b}_i \mathbf{b}_i^H]_{n,n} - p_n \right)$$

$$= \sum_{k=1}^{K} \left\{ \mathbf{b}_k^H (\sum_{i=1}^{K} \eta_i \widehat{\mathbf{h}}_i w_i w_i^H \widehat{\mathbf{h}}_i^H + \eta_i w_i w_i^H \mathbf{R}_{bi}) \mathbf{b}_k \right.$$
$$- \eta_k w_k^H \widehat{\mathbf{h}}_k^H \mathbf{b}_k - \eta_k \mathbf{b}_k^H \widehat{\mathbf{h}}_k w_k + \sigma_k^2 \eta_k w_k^H w_k + \eta_k \right\}$$
$$+ \sum_{n=1}^{N} \lambda_n \left( [\sum_{i=1}^{K} \mathbf{b}_i \mathbf{b}_i^H]_{n,n} - p_n \right)$$

$$= \sum_{k=1}^{K} \left\{ \mathbf{b}_k^H \mathbf{A} \mathbf{b}_k - \eta_k w_k^H \widehat{\mathbf{h}}_k^H \mathbf{b}_k - \eta_k \mathbf{b}_k^H \widehat{\mathbf{h}}_k w_k + \right.$$
$$\left. \sigma_k^2 \eta_k w_k^H w_k + \eta_k \right\} - \sum_{n=1}^{N} \lambda_n p_n \quad (8)$$

where $\boldsymbol{\lambda} = \text{diag}(\lambda_1, \cdots, \lambda_N)$ and $\mathbf{A} = \sum_{i=1}^{K} \eta_i \widehat{\mathbf{h}}_i w_i w_i^H \widehat{\mathbf{h}}_i^H + \eta_i w_i w_i^H \mathbf{R}_{bi} + \boldsymbol{\lambda}$. Thus, the dual function of (5) is

$$g(\boldsymbol{\lambda}) = \min_{\{\mathbf{b}_k\}_{k=1}^K} L(\boldsymbol{\lambda}, \mathbf{B}) \quad (9)$$

$$= \min_{\{\mathbf{b}_k\}_{k=1}^K} \sum_{k=1}^{K} \left\{ \mathbf{b}_k^H \mathbf{A} \mathbf{b}_k - \eta_k w_k^H \widehat{\mathbf{h}}_k^H \mathbf{b}_k - \eta_k \mathbf{b}_k^H \widehat{\mathbf{h}}_k w_k \right.$$
$$\left. + \sigma_k^2 \eta_k w_k^H w_k + \eta_k \right\} - \sum_{n=1}^{N} \lambda_n p_n$$

$$= \sum_{k=1}^{K} \left\{ \eta_k (\sigma_k^2 w_k^H w_k + 1) - \eta_k^2 w_k^H \widehat{\mathbf{h}}_k^H \mathbf{A}^{-1} \widehat{\mathbf{h}}_k w_k \right\}$$
$$- \sum_{n=1}^{N} \lambda_n p_n$$

where the third equality is obtained after substituting the optimal $\mathbf{b}_k$ of (9) which is given by

$$\mathbf{b}_k^\star = \eta_k \mathbf{A}^{-1} \widehat{\mathbf{h}}_k w_k, \; \forall k \; \Rightarrow \; \mathbf{b}_{lk}^\star = \eta_k [\mathbf{A}^{-1}]_l \widehat{\mathbf{h}}_k w_k \; \forall l, k \quad (10)$$

where $[\mathbf{A}^{-1}]_l \in C^{N_l \times N}$ is the submatrix of $\mathbf{A}^{-1}$ which is given by $[\mathbf{A}^{-1}]_l = [\mathbf{A}^{-1}]_{(F_l : F_l + N_l - 1, :)}$ with $F_l = \sum_{i=0}^{l-1} N_i + 1$ and $N_0 = 0$. As can be seen from (10), for a given $\boldsymbol{\lambda}$, the precoder of each user can be optimized independently. The optimal $\boldsymbol{\lambda}$ of (8) can be obtained by solving the dual optimization problem of (5) which is given as

$$\max_{\{\lambda_n \geq 0\}_{n=1}^N} g(\boldsymbol{\lambda}) =$$

$$\max_{\{\lambda_n \geq 0\}_{n=1}^N} \sum_{k=1}^{K} \left\{ \eta_k (\sigma_k^2 w_k^H w_k + 1) - \eta_k^2 w_k^H \widehat{\mathbf{h}}_k^H \mathbf{A}^{-1} \widehat{\mathbf{h}}_k w_k \right\}$$
$$- \sum_{n=1}^{N} \lambda_n p_n. \quad (11)$$

Considering the eigenvalue decomposition of $\widehat{\mathbf{H}} \mathbf{W} \eta^2 \mathbf{W}^H \widehat{\mathbf{H}}^H \triangleq \widetilde{\mathbf{V}} \widetilde{\boldsymbol{\Lambda}} \widetilde{\mathbf{V}}^H$ and $\widehat{\mathbf{H}} \mathbf{W} \eta \mathbf{W}^H \widehat{\mathbf{H}}^H + \sum_{i=1}^{K} \eta_i w_i w_i^H \mathbf{R}_{bi} \triangleq \bar{\mathbf{V}} \bar{\boldsymbol{\Lambda}} \bar{\mathbf{V}}^H$, problem (11) can be written as

$$\min_{\{\lambda_n \geq 0\}_{n=1}^N} \text{tr} \left\{ \mathbf{F}^H (\mathbf{R} \mathbf{R}^H + \boldsymbol{\lambda})^{-1} \mathbf{F} \right\} + \sum_{n=1}^{N} \lambda_n p_n \quad (12)$$

where $\mathbf{F} = \widetilde{\mathbf{V}} \sqrt{\widetilde{\boldsymbol{\Lambda}}}$ and $\mathbf{R} = \bar{\mathbf{V}} \sqrt{\bar{\boldsymbol{\Lambda}}}$. The above optimization problem can be cast as a semi-definite programming (SDP) problem where the global optimal solution can be found by existing convex optimization tools [18]. The worst-case computational complexity of this problem is on the order of $O((2N^2 + N)^2 (4N)^{2.5})$ [29]. However, here our aim is to obtain the optimal values of $\{\lambda_n\}_{n=1}^{N}$ distributively with less computational cost than that of the SDP method. In this regard, we present the following Lemma.

*Lemma 1*: The optimal $\{\lambda_n\}_{n=1}^{N}$ of the above optimization problem can be obtained by solving the following problem

$$\min_{\{\lambda_n, \mathbf{g}_n, \mathbf{t}_n\}_{n=1}^N} \sum_{n=1}^{N} \{\mathbf{g}_n^H \boldsymbol{\lambda}^{-1} \mathbf{g}_n + \mathbf{t}_n^H \mathbf{t}_n + \lambda_n p_n\} \triangleq \varphi$$
$$\text{s.t } \mathbf{R} \mathbf{t}_n + \mathbf{g}_n = \mathbf{f}_n, \; \forall n \quad (13)$$

where $\mathbf{f}_n$ is the $n$th column of $\mathbf{F}$.

*Proof*: By keeping $\boldsymbol{\lambda}$ constant, the Lagrangian function of (13) is given by

$$L = \sum_{n=1}^{N} \{\mathbf{g}_n^H \boldsymbol{\lambda}^{-1} \mathbf{g}_n + \mathbf{t}_n^H \mathbf{t}_n + \lambda_n p_n -$$
$$\boldsymbol{\tau}_n^H (\mathbf{R} \mathbf{t}_n + \mathbf{g}_n - \mathbf{f}_n)\} \quad (14)$$

where $\boldsymbol{\tau}_n^H$ is the Lagrangian multiplier associated with the $n$th equality constraint of (13). Differentiation of $L$ with respect to $\{\mathbf{g}_i, \mathbf{t}_i\}_{i=1}^{N}$ yield $\{\mathbf{g}_i^\star = \boldsymbol{\lambda} \boldsymbol{\tau}_i\}_{i=1}^{N}$ and $\{\mathbf{t}_i^\star = \mathbf{R}^H \boldsymbol{\tau}_i\}_{i=1}^{N}$. By substituting these $\{\mathbf{g}_i^\star, \mathbf{t}_i^\star\}_{i=1}^{N}$ in the equality constraint of (13), we get $\{\boldsymbol{\tau}_i = (\mathbf{R} \mathbf{R}^H + \boldsymbol{\lambda})^{-1} \mathbf{f}_i\}_{i=1}^{N}$. It follows

$$\mathbf{g}_i^\star = \boldsymbol{\lambda} (\mathbf{R} \mathbf{R}^H + \boldsymbol{\lambda})^{-1} \mathbf{f}_i, \; \mathbf{t}_i^\star = \mathbf{R}^H (\mathbf{R} \mathbf{R}^H + \boldsymbol{\lambda})^{-1} \mathbf{f}_i, \; \forall i \quad (15)$$

Plugging the above optimal values into the objective function of (13) yields

$$\varphi = \sum_{i=1}^{N} \{\mathbf{g}_i^H \boldsymbol{\lambda}^{-1} \mathbf{g}_i + \mathbf{t}_i^H \mathbf{t}_i + \lambda_i p_i\}$$
$$= \sum_{i=1}^{N} \{\mathbf{f}_i^H (\mathbf{R} \mathbf{R}^H + \boldsymbol{\lambda})^{-1} \mathbf{f}_i\} + \sum_{i=1}^{N} \lambda_i p_i$$
$$= \text{tr} \left\{ \mathbf{F}^H (\mathbf{R} \mathbf{R}^H + \boldsymbol{\lambda})^{-1} \mathbf{F} \right\} + \sum_{i=1}^{N} \lambda_i p_i. \quad (16)$$

The above equation is the same as the objective function of the original optimization problem (12). It follows that (12) and (13) are equivalent problems. Note that *Lemma 1* is proved by modifying the idea of matrix fractional minimization (see [18] and [20]). It can be shown that (13) is a convex optimization problem [18].

To develop distributed algorithm for (13), we reexpress $\mathbf{G} = [\mathbf{g}_1, \cdots, \mathbf{g}_N]$ as $\mathbf{G} = [\bar{\mathbf{g}}_1^H; \cdots; \bar{\mathbf{g}}_N^H]$, where $\bar{\mathbf{g}}_i^H$ is the $i$th row of $\mathbf{G}$. By doing so $\mathbf{G}^\star = [\mathbf{g}_1^\star, \cdots, \mathbf{g}_N^\star]$ of (15) can also be written as $\mathbf{G}^\star = [(\bar{\mathbf{g}}_1^\star)^H; \cdots; (\bar{\mathbf{g}}_N^\star)^H]$, where

$$\bar{\mathbf{g}}_i^\star = \lambda_i \boldsymbol{\Gamma}_i^H, \; \forall i \quad (17)$$

and $\mathbf{\Gamma}_i$ is the $i$th row of $\mathbf{\Gamma} = \mathbf{A}^{-1}\mathbf{F}$. Now, we solve (13) distributively as follows. First, keeping $\boldsymbol{\lambda}$ constant, the optimal $\overline{\mathbf{g}}_i$ can be computed independently using (17). Then, with this $\overline{\mathbf{g}}_i^\star$, $\lambda_i$ is updated by

$$\frac{\partial \varphi}{\lambda_i} = -\frac{1}{\lambda_i^2}\beta_i + p_i = 0 \Rightarrow \lambda_i^\star = \sqrt{\beta_i/p_i}, \forall i \quad (18)$$

where $\beta_i = (\overline{\mathbf{g}}_i^\star)^H \overline{\mathbf{g}}_i^\star$. As we can see from the above expression $\lambda_i^\star$ is always non-negative. Furthermore, from (17) and (18), one can observe that $\lambda_i^\star$ can be updated in parallel by using only $\overline{\mathbf{g}}_i^\star$. Thus, for our problem (12), the computation of $\{\mathbf{t}_i^\star, \mathbf{g}_i^\star\}_{i=1}^N$ is not required. To summarize, problem (12) can be solved iteratively in a distributed manner as shown in **Algorithm I**.

**Algorithm I**: Iterative algorithm to solve (12)
1) Initialization: Set $\{\lambda_n = 1\}_{n=1}^N$.
   **Repeat**
2) With the current $\boldsymbol{\lambda}$, compute $\{\overline{\mathbf{g}}_n\}_{n=1}^N$ using (17) and update $\{\lambda_n\}_{n=1}^N$ with (18).
3) Share the above $\{\lambda_n\}_{n=1}^N$ among all processors.
4) Calculate the objective function of (12).
   **Until** convergence.

As we can see, **Algorithm I** is developed to get the minimum value of the objective function of (12). Thus, this algorithm should stop iteration when the objective function of (12) is not decreasing significantly [5]. One simple approach of doing this is to stop **Algorithm I** from iteration when $\phi_i - \phi_{i+1} < \delta$, where $\phi_i$ is the objective function of (12) at the $i$th iteration and $\delta$ is the desired accuracy. For our simulation results, we have used the latter approach to declare convergence of **Algorithm I** with $\delta = 10^{-12}$.

**Convergence:** Since (13) is jointly convex in $\{\mathbf{g}_n, \mathbf{t}_n\}_{n=1}^N$ and $\{\lambda_n\}_{n=1}^N$, at each step of **Algorithm I**, the objective function of (13) is non-increasing. This implies that at each iteration of this algorithm, the objective function of (12) is also non-increasing. Moreover, it is clearly seen that the objective function of (12) is lower bounded by 0. These two facts show that **Algorithm I** is always convergent[3]. Although we are not able to prove the global optimality of **Algorithm I** analytically, in all of our simulation results, we have observed that the optimal $\boldsymbol{\lambda}$ of (12) obtained by **Algorithm I** and the SDP method are the same.

**Computational complexity:** The major computational load of **Algorithm I** arises from matrix inversion. According to [30], matrix inversion has a complexity on the order of $O(N^{2.376})$. Thus, **Algorithm I** requires $O(N^{2.376})$ per iteration. As will be shown later in Section VI, in all of our simulations, **Algorithm I** converges to an optimal solution in less than 10 iterations. This shows that the proposed distributed algorithm significantly reduces the computational load of the precoder design for $\mathcal{P}1$.

**Note:** When the $n$th power constraint of (5) is inactive, at optimality, the corresponding Lagrangian multiplier should be zero. However, when we reformulate (12) into (13), each of $\{\lambda_n\}_{n=1}^N$ is not allowed to be zero. This shows that the development of distributed algorithm for (12) with $\{\lambda_n \geq 0\}_{n=1}^N$ is an open problem.

Using $\{\lambda_n\}_{n=1}^N$ of **Algorithm I**, the optimal $\{\mathbf{b}_{lk}\}_{k=1}^K, \forall l$ of (5) can be computed by (10), and with these $\{\mathbf{b}_{lk}\}_{k=1}^K, \forall l$, the decoder of each user can be computed by using the MAMSE receiver approach (4). In summary, we solve $\mathcal{P}1$ (5) distributively as shown in **Algorithm II**.

**Algorithm II**: Distributed algorithm for problem $\mathcal{P}1$ (5)
Initialization: Set $\mathbf{B} = \hat{\mathbf{H}}$ and normalize the rows of $\mathbf{B}$ such that the power constraint of each antenna is satisfied with equality. Then, initialize $\{w_k\}_{k=1}^K$ by MAMSE receiver (4). Set the maximum number of iterations $i_{max}$.
**Repeat**
1) Compute the optimal $\{\lambda_n\}_{n=1}^N$ with **Algorithm I**.
2) Solve for $\{\mathbf{b}_k\}_{k=1}^K$ using (10).
3) Update the MAMSE receivers $\{w_k\}_{k=1}^K$ with (4).
4) Calculate the objective function of (5).
   **Until** convergence.

In our simulation, we declare the convergence of this algorithm when $\bar{\xi}_i - \bar{\xi}_{i+1} < 10^{-6}$, where $\bar{\xi}_i$ is the achieved weighted sum AMSE at the $i$th iteration of **Algorithm II**.

**Convergence:** It can be shown that at each iteration of **Algorithm II**, the objective function of (5) is non-increasing. Since we are interested to get any local optimal $\{\mathbf{b}_k, w_k\}_{k=1}^K$ that yields the local minimum $\bar{\xi}$, the convergence analysis of **Algorithm II** with respect to the optimization variables $\{\mathbf{b}_k, w_k\}_{k=1}^K$ is not required.

**Implementation of Algorithm II:** For the implementation of this distributed algorithm, for simplicity, it is assumed that $L = N = K$, and $\mathcal{P}1$ is solved in a central controller which has $K$ parallel processors. We can implement **Algorithm II** distributively as follows.

Initialization: Each processor sets $\{w_k\}_{k=1}^K$ as in **Algorithm II** and $\{\lambda_n = 1\}_{n=1}^N$.
1) With the current $\{\lambda_n\}_{n=1}^N$ and $\{w_k\}_{k=1}^K$, the $n$th processor computes $\overline{\mathbf{g}}_n$ using (17) and updates its $\lambda_n$ by (18), $\forall n$. Then, $\{\lambda_n\}_{n=1}^N$ are shared among all processors. These two steps are repeated until $\{\lambda_n\}_{n=1}^N$ are found to be optimal.
2) Using $\{\lambda_n\}_{n=1}^N$ of step 1, the $k$th processor computes the optimal $\mathbf{b}_k$ by (10), $\forall k$ and $\{\mathbf{b}_k\}_{k=1}^K$ are shared among all processors. Again, using these precoders, the $k$th processor computes $w_k$ with (4), $\forall k$ and $\{w_k\}_{k=1}^K$ are shared to processors.
3) Steps (1) and (2) are repeated until **Algorithm II** is convergent.
4) The controller finally sends the optimal precoders and decoders to the corresponding BSs and MSs, respectively.

Note that in some scenario we might be obliged to design the precoders and decoders of all users without a central controller. In this case, one can apply the above implementation approach just by replacing the role of processors with that of BSs.

---
[3]Note that since the aim of problem (12) is to get any $\{\lambda_n\}_{n=1}^N$ which achieves the smallest objective function of (12), we believe that the convergence analysis of **Algorithm I** with respect to the optimization variables $\{\lambda_n\}_{n=1}^N$ is not required.

## B. Robust power minimization problem ($\mathcal{P}2$)

The robust power minimization constrained with the MSE target of each user and the power of each BS antenna problem is formulated as

$$\min_{\{\mathbf{b}_k, w_k\}_{k=1}^K} \sum_{k=1}^K \mathbf{b}_k^H \mathbf{b}_k,$$

$$\text{s.t } [\sum_{k=1}^K \mathbf{b}_k \mathbf{b}_k^H]_{n,n} \leq p_n, \ \bar{\xi}_k \leq \varepsilon_k, \ 0 < \varepsilon_k < 1, \ \forall k \quad (19)$$

where $\varepsilon_k$ is the $k$th user AMSE target.

*1) Centralized precoder design for $\mathcal{P}2$:* For fixed receivers $\{w_k\}_{k=1}^K$, using (3) and applying the same technique as in (7), the above problem can be equivalently expressed as [17]

$$\min_{\mathbf{B}, \tilde{\chi}} \tilde{\chi} \quad (20)$$

$$\text{s.t } \|\text{vec}(\mathbf{B})\|_2 \leq \tilde{\chi}, \ \|\tilde{\mathbf{b}}_n\|_2 \leq \sqrt{p_n}, \ \forall n$$

$$\|[(\mathbf{B}^H \widehat{\mathbf{h}}_k w_k - \boldsymbol{\theta}_k); \text{vec}(\sqrt{\mathbf{R}_{bk}} \mathbf{B} w_k)]\|_2 \leq$$

$$\sqrt{\varepsilon_k - \sigma_k^2 w_k^H w_k}, \ \forall k$$

where $\boldsymbol{\theta}_k$ is a column vector of size $K$ with the $k$th element equal to 1 and all the other elements equal to 0. The above problem is a SOCP for which the global optimal solution can be obtained with $O(\sqrt{(K+N+1)}(2KN+1)^2(2K^2+2NK^2+4NK))$ computational load [29].

*2) Distributed precoder design for $\mathcal{P}2$:* For fixed $\{w_k\}_{k=1}^K$, like in $\mathcal{P}1$, we utilize the Lagrangian dual decomposition method to solve $\mathcal{P}2$ distributively. The Lagrangian function of (19) is given as

$$L(\boldsymbol{\lambda}, \tilde{\boldsymbol{\nu}}, \mathbf{B}) = \sum_{k=1}^K \{\mathbf{b}_k^H (\mathbf{I}_N + \sum_{i=1}^K \nu_i \widehat{\mathbf{h}}_i w_i w_i^H \widehat{\mathbf{h}}_i^H + \nu_i w_i w_i^H \mathbf{R}_{bi}$$

$$+ \boldsymbol{\lambda}) \mathbf{b}_k - 2\Re\{\nu_k w_k^H \widehat{\mathbf{h}}_k^H \mathbf{b}_k\} + \sigma_k^2 \nu_k w_k^H w_k +$$

$$\nu_k(1-\varepsilon_k)\} - \sum_{n=1}^N \lambda_n p_n \quad (21)$$

where $\boldsymbol{\lambda}$ and $\tilde{\boldsymbol{\nu}} = \text{diag}(\nu_1, \cdots, \nu_K)$ are the Lagrangian multipliers for the first and second constraint sets of (19), respectively. Thus, the dual function of (19) is computed by

$$g(\boldsymbol{\lambda}, \tilde{\boldsymbol{\nu}}) = \min_{\{\mathbf{b}_k\}_{k=1}^K} L(\boldsymbol{\lambda}, \tilde{\boldsymbol{\nu}}, \mathbf{B})$$

$$= \sum_{k=1}^K \nu_k \alpha_k - \sum_{k=1}^K \nu_k w_k^H \widehat{\mathbf{h}}_k^H \widetilde{\mathbf{A}}^{-1} \widehat{\mathbf{h}}_k w_k \nu_k - \sum_{n=1}^N \lambda_n p_n \quad (22)$$

where $\alpha_k = 1 + \sigma_k^2 w_k^H w_k - \varepsilon_k$, $\widetilde{\mathbf{A}} = \mathbf{I}_N + \sum_{i=1}^K \nu_i \widehat{\mathbf{h}}_i w_i w_i^H \widehat{\mathbf{h}}_i^H + \nu_i w_i w_i^H \mathbf{R}_{bi} + \boldsymbol{\lambda}$ and the second equality is obtained after substituting the optimum $\mathbf{b}_k$ which is given by

$$\mathbf{b}_k^\star = \nu_k \widetilde{\mathbf{A}}^{-1} \widehat{\mathbf{h}}_k w_k, \ \forall k. \quad (23)$$

From (4) and (23) it can be clearly seen that if (19) is feasible, $\{\nu_k > 0\}_{k=1}^K$ must be satisfied. Furthermore, for given $\{\lambda_n\}_{n=1}^N$ and $\{\nu_k\}_{k=1}^K$, the precoder of each user can be optimized distributively by using (23). The optimal $\{\lambda_n\}_{n=1}^N$ and $\{\nu_k\}_{k=1}^K$ of (21) can be obtained by solving the dual problem of (19) which can be expressed as

$$\max_{\{\lambda_n \geq 0\}_{n=1}^N, \{\nu_k > 0\}_{k=1}^K} g(\boldsymbol{\lambda}, \tilde{\boldsymbol{\nu}}) =$$

$$\max_{\{\lambda_n \geq 0\}_{n=1}^N, \{\nu_k > 0\}_{k=1}^K} \sum_{k=1}^K \nu_k \alpha_k - w_k^H \widehat{\mathbf{h}}_k^H \widetilde{\mathbf{A}}^{-1} \widehat{\mathbf{h}}_k w_k \nu_k^2$$

$$- \sum_{n=1}^N \lambda_n p_n. \quad (24)$$

It can be shown that the above problem can be formulated as SDP [18]. Thus, (24) can be solved by using convex optimization tools with complexity on the order of $O(\sqrt{K(N+1)}(N+2K)^2 K(N+1)^2)$ [29]. However, our interest is to obtain the optimal values of $\{\lambda_n\}_{n=1}^N$ and $\{\nu_k\}_{k=1}^K$ of the above problem distributively. The above problem can be rewritten as

$$\min_{\{\lambda_n \geq 0\}_{n=1}^N, \{\nu_k > 0\}_{k=1}^K} \sum_{k=1}^K \nu_k w_k^H \widehat{\mathbf{h}}_k^H (\bar{\mathbf{W}} \boldsymbol{\Upsilon} \bar{\mathbf{W}}^H + \mathbf{I})^{-1} \widehat{\mathbf{h}}_k w_k \nu_k$$

$$+ \sum_{n=1}^N \lambda_n + \sum_{k=1}^K \nu_k \alpha_k p_n \quad (25)$$

where $\boldsymbol{\Upsilon} = \text{blkdiag}(\boldsymbol{\nu}, \boldsymbol{\lambda})$, $\bar{\mathbf{W}} = [\widetilde{\mathbf{W}} \ \mathbf{I}_N]$, $\boldsymbol{\nu} = \text{blkdiag}(\nu_1 \mathbf{I}_N, \cdots, \nu_K \mathbf{I}_N)$, $\widetilde{\mathbf{W}} = [\widetilde{\mathbf{W}}_{b1}, \cdots, \widetilde{\mathbf{W}}_{bK}]$ and $\widetilde{\mathbf{W}}_{bk} = (w_k w_k^H (\widehat{\mathbf{h}}_k \widehat{\mathbf{h}}_k^H + \mathbf{R}_{bk}))^{1/2}$. Now, by applying matrix fractional minimization of [18] on the first sum terms of the above problem, we can reformulate (25) as (see page 198 of [18])

$$\min_{\{\bar{\mathbf{g}}_k, \bar{\mathbf{t}}_k, \nu_k\}_{k=1}^K, \{\lambda_n\}_{n=1}^N} \sum_{k=1}^K \bar{\mathbf{g}}_k^H \boldsymbol{\Upsilon}^{-1} \bar{\mathbf{g}}_k + \bar{\mathbf{t}}_k^H \bar{\mathbf{t}}_k + \sum_{n=1}^N \lambda_n p_n$$

$$- \sum_{k=1}^K \alpha_k \nu_k, \quad \text{s.t } \bar{\mathbf{t}}_k + \bar{\mathbf{W}} \bar{\mathbf{g}}_k = \nu_k \widehat{\mathbf{h}}_k w_k, \ \forall k. \quad (26)$$

From the equality constraint of (26), we get $\bar{\mathbf{t}}_k = \nu_k \widehat{\mathbf{h}}_k w_k - \bar{\mathbf{W}} \bar{\mathbf{g}}_k$. By substituting this $\bar{\mathbf{t}}_k$ into the objective function of the above problem, (26) can be rewritten as

$$\min_{\{\bar{\mathbf{g}}_k, \nu_k\}_{k=1}^K, \{\lambda_n\}_{n=1}^N} \sum_{k=1}^K \bar{\mathbf{g}}_k^H \boldsymbol{\Upsilon}^{-1} \bar{\mathbf{g}}_k + (\nu_k \widehat{\mathbf{h}}_k w_k - \bar{\mathbf{W}} \bar{\mathbf{g}}_k)^H.$$

$$(\nu_k \widehat{\mathbf{h}}_k w_k - \bar{\mathbf{W}} \bar{\mathbf{g}}_k) + \sum_{n=1}^N \lambda_n p_n - \sum_{k=1}^K \alpha_k \nu_k. \quad (27)$$

For fixed $\{\nu_k\}_{k=1}^K$ and $\{\lambda_n\}_{n=1}^N$, the optimal $\bar{\mathbf{g}}_k$ of problem (27) is given by

$$\bar{\mathbf{g}}_k^\star = \nu_k (\boldsymbol{\Upsilon}^{-1} + \bar{\mathbf{W}}^H \bar{\mathbf{W}})^{-1} \bar{\mathbf{W}}^H \widehat{\mathbf{h}}_k w_k$$

$$= \nu_k \boldsymbol{\Upsilon} \bar{\mathbf{W}}^H \widetilde{\mathbf{A}}^{-1} \widehat{\mathbf{h}}_k w_k, \forall k \quad (28)$$

where the second equality is obtained by employing matrix inversion Lemma [28]. To develop distributed algorithm for the above problem, we introduce the following variables: $\{\widetilde{\mathbf{G}}_k^\star \in \mathcal{C}^{N \times K} \text{ as the } \bar{\mathbf{G}}_{[N(k-1)+1:Nk,:]}^\star\}_{k=1}^K$ submatrix of $\bar{\mathbf{G}}^\star = [\bar{\mathbf{g}}_1^\star, \cdots, \bar{\mathbf{g}}_K^\star]$ and $\{(\bar{\mathbf{u}}_n^\star)^H \text{ as the } KN+n\}_{n=1}^N$ row of $\bar{\bar{\mathbf{G}}}^\star$. For the given $\{\nu_k\}_{k=1}^K$ and $\{\lambda_n\}_{n=1}^N$, $\tilde{\mathbf{G}}_k^\star$ and $(\bar{\mathbf{u}}_n^\star)^H$ can also be computed as

$$\tilde{\mathbf{G}}_k^\star = \nu_k \tilde{\mathbf{W}}_{bk}^H \tilde{\mathbf{\Gamma}}, \ \forall k, \quad \bar{\mathbf{u}}_n^\star = \lambda_n \tilde{\mathbf{\Gamma}}_n^H, \ \forall n \quad (29)$$

where $\tilde{\mathbf{\Gamma}} = \tilde{\mathbf{A}}^{-1} \widehat{\mathbf{H}} \mathbf{W} \tilde{\boldsymbol{\nu}}$ and $\tilde{\mathbf{\Gamma}}_n$ is the $n$th row of $\tilde{\mathbf{\Gamma}}$.

Now, we solve problem (27) distributively as follows. First, for fixed $\{\nu_k\}_{k=1}^K$ and $\{\lambda_n\}_{n=1}^N$, the optimal $\bar{\bar{\mathbf{g}}}_k$ of (27) and the introduced variables ($\tilde{\mathbf{G}}_k^\star, \bar{\mathbf{u}}_n^\star$) are computed using (28) and (29), respectively. Then, using these $\bar{\bar{\mathbf{g}}}_k^\star$, $\tilde{\mathbf{G}}_k^\star$ and $\bar{\mathbf{u}}_n^\star$, $\nu_k$ and $\lambda_n$ are updated independently and distributively by

$$\nu_k^\star = \min_{\nu_k > 0} \ \nu_k^2 \rho_{k1} - \nu_k \rho_{k2} + \frac{\rho_{k3}}{\nu_k} \quad (30)$$

$$\lambda_n^\star = \sqrt{\frac{\rho_{n0}}{p_n}}, \ \forall n \quad (31)$$

where $\rho_{k1} = \widehat{\mathbf{h}}_k^H w_k^H w_k \widehat{\mathbf{h}}_k$, $\rho_{k2} = 2\Re\{w_k^H \widehat{\mathbf{h}}_k^H \bar{\mathbf{W}} \bar{\bar{\mathbf{g}}}_k^\star\} + \alpha_k$, $\rho_{k3} = \text{tr}\{\tilde{\mathbf{G}}_k^\star (\tilde{\mathbf{G}}_k^\star)^H\}$ and $\rho_{n0} = (\bar{\mathbf{u}}_n^\star)^H \bar{\mathbf{u}}_n^\star$. If $\rho_{k1} \neq 0$, by applying first order derivative, it can be shown that (30) has exactly one real solution which is given by [31]

$$\nu_k^\star = \frac{1}{6\rho_{k1}} \left[\rho_{k2} + \mu_{k1} + \mu_{k2}\right], \ \forall k \quad (32)$$

where $\mu_{k1} = \sqrt[3]{\frac{1}{2}\left[2\rho_{k2}^3 + c_k - \zeta_k\right]}$, $\mu_{k2} = \sqrt[3]{\frac{1}{2}\left[2\rho_{k2}^3 + c_k + \zeta_k\right]}$, $c_k = 108\rho_{k1}^2 \rho_{k3}$ and $\zeta_k = \sqrt{c_k(c_k + 4\rho_{k2}^3)}$. One can easily see that $\rho_{k1} \geq 0$, $\mu_{k1} \geq 0$, $\mu_{k2} \geq 0$ and $\rho_{k2} > 0$. Moreover, when $\rho_{k1} > 0$, $\nu_k^\star$ of (32) is always positive. To summarize, (25) can be solved distributively in an iterative manner as in **Algorithm III**.

**Algorithm III**: Distributed algorithm to solve (25)
1) Initialization: Set $\{\lambda_n = 1\}_{n=1}^N$ and $\{\nu_k = 1\}_{k=1}^K$.
   **Repeat**
2) With the current $\{\lambda_n\}_{n=1}^N$ and $\{\nu_k\}_{k=1}^K$, compute $\bar{\bar{\mathbf{g}}}_k, \tilde{\mathbf{G}}_k$ and $\bar{\mathbf{u}}_n$ using (28) and (29), $\forall k, n$. Then, update $\lambda_n$ and $\nu_k$ by (31) and (32), respectively, $\forall n, k$.
3) Share the above $\{\lambda_n\}_{n=1}^N$ and $\{\nu_k\}_{k=1}^K$ among all processors.
4) Calculate the objective function of (25).
   **Until** convergence.

**Feasibility study for $\mathcal{P}2$:** The problem (19) is infeasible, if there exists at least one MS with either $\rho_{k1} = 0$ or $\nu_k^\star \gg 0$. This can be justified as follows. For the former case (i.e., $\exists k$ such that $\rho_{k1} = 0$), one can easily verify that (19) is infeasible. For the latter case (i.e., when $\{\rho_{k1} > 0\}_{k=1}^K$ and $\exists k$ such that $\nu_k^\star \gg 0$), although $\nu_k^\star$ are not permitted to be $\infty$, $\nu_k^\star$ can be arbitrarily very large number. And when $\nu_k^\star$ is large, one can use (23) to show that the $k$th MS needs additional power at least in one of the BS antennas to satisfy its AMSE target. This case corresponds to the scenario where (25) is an unbounded problem. During the iterative stages of **Algorithm III**, when (19) is infeasible, we have observed from the simulation results that $\{\nu_k\}_{k=1}^K$ of (32) increases rapidly for at least one MS and the latter algorithm never converges to a point.

**Convergence:** If $\mathcal{P}2$ is feasible, it can be shown that **Algorithm III** is guaranteed to converge. However, we are not able to show the global optimality of **Algorithm III** analytically. Nonetheless, in all simulation results we observe that the optimal $\boldsymbol{\lambda}$ and $\tilde{\boldsymbol{\nu}}$ of (25) obtained by the latter algorithm and SDP method are the same.

**Computational complexity:** As can be seen from (28) and (29), in the proposed distributed algorithm, the main computational load comes from the computation of $\tilde{\mathbf{A}}^{-1}$ which can also be computed efficiently with $O(N^{2.376})$ [30]. Moreover, in all of our simulation results, we have observed that **Algorithm III** converges to an optimal solution within few iterations.

Once we get the optimal $\{\lambda_n\}_{n=1}^N$ and $\{\nu_k\}_{k=1}^K$, like in $\mathcal{P}1$, the precoders and decoders of each user can be optimized using (23) and (4), respectively. It follows that $\mathcal{P}2$ (19) can be solved distributively like in **Algorithm II** of $\mathcal{P}1$.

## V. EXTENSION TO MIMO COORDINATED BASE STATION SYSTEMS

For multiuser MIMO coordinated BS systems, the solution approaches of Section IV can be applied to solve the following problems. (1) The robust minimization of symbol (user) wise weighted sum MSE with per BS antenna power constraint problem. (2) The robust minimization of the total sum power of all BSs with per BS antenna power and per symbol (user) MSE target constraints. In this section, we examine the robust symbol wise weighted sum MSE minimization with a per BS antenna power constraint problem for the multiuser MIMO coordinated BS systems ($\tilde{\mathcal{P}}1$) only[4]. For the MIMO coordinated BS systems, the channel estimation technique of Section III can be utilized. Upon doing so, the true channel between the $l$th BS and the $k$th user $\mathbf{H}_{lk}^H$ and its MMSE estimate $\widehat{\mathbf{H}}_{lk}^H$ are related by [11]

$$\mathbf{H}_{lk}^H = \widehat{\mathbf{H}}_{lk}^H + \mathbf{E}_{wlk}^H \tilde{\mathbf{R}}_{blk}^{1/2} = \widehat{\mathbf{H}}_{lk}^H + \mathbf{E}_{lk}^H \quad (33)$$

where $\mathbf{E}_{lk}^H \in \mathcal{C}^{M_k \times N_l}$ and $M_k$ are the estimation error matrix and number of antennas of the $k$th MS, respectively, and the entries of $\mathbf{E}_{wlk}^H$ are i.i.d with $\mathcal{CN}(0, \sigma_{elk}^2)$. Like in (3), the AMSE of the $k$th MS $i$th symbol ($\bar{\xi}_{ki}$) can be expressed as

$$\bar{\xi}_{ki} = 1 + \mathbf{w}_{ki}^H (\widehat{\mathbf{H}}_k^H \sum_{j=1}^K \sum_{m=1}^{S_j} \mathbf{b}_{jm} \mathbf{b}_{jm}^H \widehat{\mathbf{H}}_k +$$
$$\text{tr}\{\mathbf{R}_{bk} \sum_{j=1}^K \sum_{m=1}^{S_j} \mathbf{b}_{jm} \mathbf{b}_{jm}^H\} \mathbf{I}_{M_k} + \sigma_k^2 \mathbf{I}_{M_k}) \mathbf{w}_{ki} -$$
$$\mathbf{b}_{ki}^H \widehat{\mathbf{H}}_k \mathbf{w}_{ki} - \mathbf{w}_{ki}^H \widehat{\mathbf{H}}_k^H \mathbf{b}_{ki} \quad (34)$$

where $\mathbf{b}_{ki} \in \mathcal{C}^{N \times 1}$ and $\mathbf{w}_{ki} \in \mathcal{C}^{M_k \times 1}$ are the precoder and decoder vectors of the $k$th MS $i$th symbol, respectively, and $\widehat{\mathbf{H}}_k^H = [\widehat{\mathbf{H}}_{1k}^H, \cdots, \widehat{\mathbf{H}}_{Lk}^H] \in \mathcal{C}^{M_k \times N}$. Using (33) and (34), we can formulate $\tilde{\mathcal{P}}1$ as

$$\min_{\{\mathbf{w}_{ki}, \mathbf{b}_{ki}\}_{k=1}^K} \sum_{k=1}^K \sum_{i=1}^{S_k} \eta_{ki} \bar{\xi}_{ki},$$
$$\text{s.t} \left[\sum_{k=1}^K \sum_{i=1}^{S_k} \mathbf{b}_{ki} \mathbf{b}_{ki}^H\right]_{n,n} \leq p_n, \ \forall n \quad (35)$$

---
[4]Note that all the other problems of this section can be examined like in $\tilde{\mathcal{P}}1$.

where $\eta_{ki}$ is the AMSE weighting factor of the $k$th MS $i$th symbol. For given precoder vectors $\{\mathbf{b}_{ki}, \forall i\}_{k=1}^K$, the receivers $\{\mathbf{w}_{ki}, \forall i\}_{k=1}^K$ of the above problem can be optimized by the following MAMSE approach

$$\mathbf{w}_{ki} = \left(\widehat{\mathbf{H}}_k^H \sum_{j=1}^K \sum_{m=1}^{S_j} \mathbf{b}_{jm}\mathbf{b}_{jm}^H \widehat{\mathbf{H}}_k + \mathrm{tr}\{\mathbf{R}_{bk}\sum_{j=1}^K\sum_{m=1}^{S_j} \mathbf{b}_{jm}\mathbf{b}_{jm}^H\}\mathbf{I}_{M_k} + \sigma_k^2 \mathbf{I}_{M_k}\right)^{-1} \widehat{\mathbf{H}}_k^H \mathbf{b}_{ki}, \forall k,i. \quad (36)$$

Next, for fixed $\{\mathbf{w}_{ki}, \forall i\}_{k=1}^K$, we summarize the centralized and distributed precoder design algorithms of $\tilde{\mathcal{P}}1$.

### A. Centralized precoder design of $\tilde{\mathcal{P}}1$

By employing (34), $\sum_{k=1}^K \sum_{i=1}^{S_k} \bar{\xi}_{ki}$ can be written as a quadratic expression like that of (6). This shows that the precoder design problem of (35) can be formulated as SOCP for which the global optimal solution can be obtained using convex optimization tools (see also [17]).

### B. Distributed precoder design of $\tilde{\mathcal{P}}1$

Here like in $\mathcal{P}1$, the precoder design problem of $\tilde{\mathcal{P}}1$ can be solved distributively by applying the Lagrangian dual decomposition and modified matrix fractional minimization approaches. After some straightforward steps, the Lagrangian function associated with $\tilde{\mathcal{P}}1$ can be expressed as

$$L(\boldsymbol{\lambda}, \{\mathbf{b}_{ki}, \forall i\}_{k=1}^K) = \sum_{k=1}^K \sum_{i=1}^{S_k} \bigg\{\mathbf{b}_{ki}^H \widetilde{\widetilde{\mathbf{A}}} \mathbf{b}_{ki} - \eta_{ki} \mathbf{w}_{ki}^H \widehat{\mathbf{H}}_k^H \mathbf{b}_{ki} - \eta_{ki}\mathbf{b}_{ki}^H \widehat{\mathbf{H}}_{ki}\mathbf{w}_{ki} + \sigma_k^2 \eta_{ki}\mathbf{w}_{ki}^H\mathbf{w}_{ki} + \eta_{ki}\bigg\} - \sum_{n=1}^N \lambda_n p_n \quad (37)$$

where $\boldsymbol{\lambda} = \mathrm{diag}(\lambda_1,\cdots,\lambda_N)$ are the Lagrangian multipliers corresponding to the constraint sets of (35) and $\widetilde{\widetilde{\mathbf{A}}} = \sum_{j=1}^K \sum_{m=1}^{S_j} \eta_{jm}\widehat{\mathbf{H}}_j \mathbf{w}_{jm}\mathbf{w}_{jm}^H \widehat{\mathbf{H}}_j^H + \mathrm{tr}\{\mathbf{w}_{jm}^H\mathbf{w}_{jm}\}\mathbf{R}_{bj} + \boldsymbol{\lambda}$. By employing the above expression and after some mathematical manipulations, the dual problem of (35) can be formulated as

$$\max_{\{\lambda_n \geq 0\}_{n=1}^N} \min_{\{\mathbf{b}_{ki},\forall i\}_{k=1}^K} L(\boldsymbol{\lambda},\{\mathbf{b}_{ki},\forall i\}_{k=1}^K) = \max_{\{\lambda_n\geq 0\}_{n=1}^N} \sum_{k=1}^K \sum_{i=1}^{S_k} \bigg\{\eta_{ki}(\sigma_k^2 \mathbf{w}_{ki}^H\mathbf{w}_{ki} + 1) - \eta_{ki}^2 \mathbf{w}_{ki}^H \widehat{\mathbf{H}}_k^H \widetilde{\widetilde{\mathbf{A}}}^{-1} \widehat{\mathbf{H}}_k\mathbf{w}_{ki}\bigg\} - \sum_{n=1}^N \lambda_n p_n. \quad (38)$$

This problem has exactly the same structure at that of (11). Thus, with the help of *Lemma 1*, we can develop distributed algorithm to solve the above problem. Consequently, $\tilde{\mathcal{P}}1$ can be solved distributively like that of $\mathcal{P}1$.

## VI. SIMULATION RESULTS

In this section, we present the simulation results for problems $\mathcal{P}1$ and $\mathcal{P}2$. The spacial antenna correlation matrix between the $l$th BS and the $k$th user $\tilde{\mathbf{R}}_{blk}$ is taken from a widely used exponential correlation model as $\{\tilde{\mathbf{R}}_{blk} = \rho_{blk}^{|i-j|}\}_{k=1}^K, \forall l$, where $0 \leq \rho_{blk} < 1$ and $1 \leq i(j) \leq N_l$, and $\{\sigma_{e1k}^2 = \sigma_{e2k}^2 =,\cdots,= \sigma_{eLk}^2 = \sigma_{ek}^2\}_{k=1}^K$. We have used exponential correlation model because of the following two reasons. First, exponential correlation model is physically reasonable in a way that the correlation between two transmit antennas decreases as the distance between them increases [32]. Second, this model is a widely used antenna correlation model for an urban area communications [26].

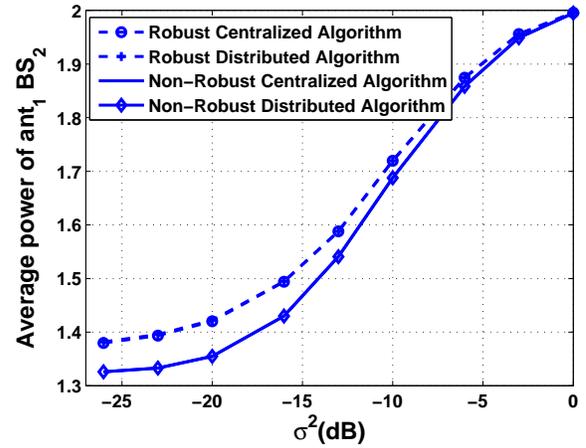

Fig. 2. The average power utilized by the first antenna of $BS_2$ for the robust centralized, robust distributed, non-robust centralized [17] and non-robust distributed [20] designs when $\rho_{b11} = 0.25, \rho_{b12} = 0.5, \rho_{b13} = 0.2, \rho_{b14} = 0.4, \rho_{b21} = 0.6, \rho_{b22} = 0.1, \rho_{b23} = 0.8$ and $\rho_{b24} = 0.15$.

### A. Simulation results for $\mathcal{P}1$

In this subsection, we consider a system with $L = 2$ BSs where each BS has 2 antennas and $K = 4$ MSs. We use $\rho_{b11} = 0.25, \rho_{b12} = 0.5, \rho_{b13} = 0.2, \rho_{b14} = 0.4, \rho_{b21} = 0.6, \rho_{b22} = 0.1, \rho_{b23} = 0.8$ and $\rho_{b24} = 0.15$, and $\sigma_{e1}^2 = 0.01, \sigma_{e2}^2 = 0.02, \sigma_{e3}^2 = 0.03$ and $\sigma_{e4}^2 = 0.04$. It is assumed that $\{\sigma_k^2 = \sigma^2\}_{k=1}^K$, $\{p_n = 2\}_{n=1}^4$ and $\{\eta_k = 1\}_{k=1}^K$. All simulation results of this subsection are averaged over 100 randomly chosen channel realizations.

The optimal transmit power of the first antenna of $BS_2$ as a function of the noise power is plotted in Fig. 2. This figure shows that the power utilized by the centralized algorithms of the robust and non-robust/naive designs are the same as that of the distributed algorithms of the robust and non-robust designs, respectively. The non-robust/naive design refers to the design in which the estimated channel is considered as perfect [17], [20]. The latter figure also shows that all antennas do not necessarily utilize their full powers to minimize the total sum AMSE of the system in the robust and non-robust designs[5]. Furthermore, the power utilized by the first antenna of $BS_2$ for both of these designs are not necessarily the same. Although the robust and non-robust designs do not use the same power at each antenna, we have noticed at all SNR values that the average total sum power of all antennas of the latter designs are almost the same. In the sequel, we compare the performance of the robust centralized, robust distributed, non-robust centralized and non-robust distributed algorithms in terms of sum AMSE. For this purpose, we define the signal-to-noise ratio (SNR) as $P_{\text{sum}}/\sigma^2$, where $P_{\text{sum}}$ is the total sum power utilized by all BS antennas of the robust distributed algorithm and $\sigma^2$ is the noise variance. The SNR is controlled by varying $\sigma^2$.

---
[5]This behavior has also been observed in [17] and [20] where the sum MSE minimization with per BS antenna power constraint problem is examined.

We first compare the performance of the robust centralized and robust distributed algorithms in terms of sum AMSE. Fig. 3 shows that the robust centralized and distributed algorithms achieve the same sum AMSE. Next, we compare the performance of the robust and non-robust/naive designs. In [20], we have shown that [17] and [20] achieve the same sum MSE. Thus, it is sufficient to compare the performance of our robust designs with the non-robust design of [20]. As can be seen from Fig. 3, the proposed robust designs have better performance than that of the non-robust design in [20] and this improvement is better at high SNR regions. As can be seen from the second equality of (3), when SNR is high (i.e., when $\sigma_k^2 = \sigma^2$ is low), $\sigma^2$ is negligible compared to $\sum_{l=1}^{L} \sigma_{elk}^2 \text{tr}\{\tilde{\mathbf{R}}_{blk}^{1/2} \mathbf{B}_l \mathbf{B}_l^H \tilde{\mathbf{R}}_{blk}^{1/2}\}$ (the term due to channel estimation error). Thus, in the high SNR region, since the non-robust design does not take into account the effect of $\sum_{l=1}^{L} \sigma_{elk}^2 \text{tr}\{\tilde{\mathbf{R}}_{blk}^{1/2} \mathbf{B}_l \mathbf{B}_l^H \tilde{\mathbf{R}}_{blk}^{1/2}\}$ which is the dominant term, the sum AMSE of this design increases significantly. These discussions help us to understand why the performance of non-robust design worsens in the high SNR region. From this explanation, we can imply that the sum AMSE gap between robust and non-robust designs increases as the SNR increases.

To see the effect of antenna correlation matrix on the achievable sum AMSE of robust and non-robust designs, we change the latter $\{\rho_{blk}\}_{k=1}^{K}, \forall l$ to $\rho_{b11} = 0.35, \rho_{b12} = 0.5, \rho_{b13} = 0.3, \rho_{b14} = 0.4, \rho_{b21} = 0.6, \rho_{b22} = 0.2, \rho_{b23} = 0.8$ and $\rho_{b24} = 0.25$. For this setting, we plot the total sum AMSE of our robust design and the non-robust design of [20] in Fig. 4. By examining this figure and Fig. 3, one can see that the sum AMSE increases as the antenna correlation coefficient increases for both the robust and non-robust designs. This is because when $\{\rho_{blk}, \exists l, k\}$ increases, the number of symbols with low channel gain increases (this can be easily seen from the eigenvalue decomposition of $\mathbf{R}_{bk}$). Consequently, for a given SNR value, the total sum AMSE also increases. The above scenario gracefully fits to that of [12] where the robust weighted sum MSE minimization with a total BS power constraint problem is examined for the conventional downlink MIMO systems.

### B. Simulation results for $\mathcal{P}2$

For the simulation result of $\mathcal{P}2$, we also consider a system with $L = 2$ BSs where each BS has 2 antennas and $K = 4$ MSs. We use $\rho_{b11} = 0.25, \rho_{b12} = 0.5, \rho_{b13} = 0.2, \rho_{b14} = 0.4, \rho_{b21} = 0.6, \rho_{b22} = 0.1, \rho_{b23} = 0.8$ and $\rho_{b24} = 0.15$, and

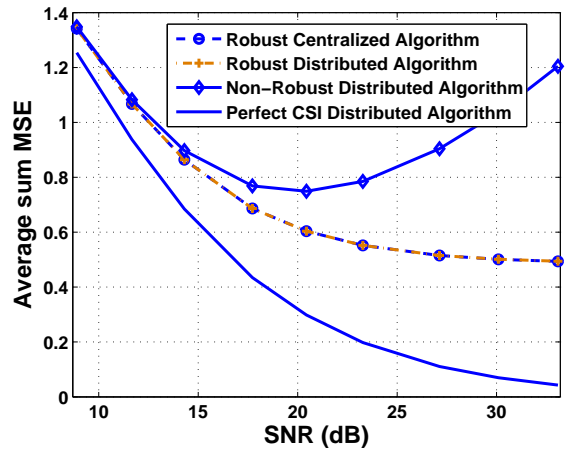

Fig. 3. Comparison of the robust centralized design, distributed design and the non-robust design of [20] when $\rho_{b11} = 0.25, \rho_{b12} = 0.5, \rho_{b13} = 0.2, \rho_{b14} = 0.4, \rho_{b21} = 0.6, \rho_{b22} = 0.1, \rho_{b23} = 0.8$ and $\rho_{b24} = 0.15$.

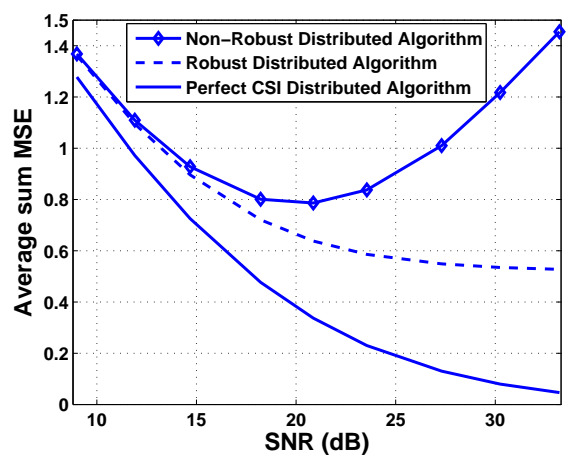

Fig. 4. Comparison of the robust design and non-robust design of [20] when $\rho_{b11} = 0.35, \rho_{b12} = 0.5, \rho_{b13} = 0.3, \rho_{b14} = 0.4, \rho_{b21} = 0.6, \rho_{b22} = 0.2, \rho_{b23} = 0.8$ and $\rho_{b24} = 0.25$.

$\sigma_{e1}^2 = 0.01$, $\sigma_{e2}^2 = 0.02$, $\sigma_{e3}^2 = 0.03$ and $\sigma_{e4}^2 = 0.04$. It is assumed that $\{\sigma_k^2 = \sigma^2\}_{k=1}^{K}$, $\{p_n = 15\}_{n=1}^{4}$ and the AMSE target of each user is set to $\{\varepsilon_k = 0.2\}_{k=1}^{K}$[6]. For better explanation of the simulation results of this subsection, we use the following channel estimate obtained from the above settings. where the rows of $\widehat{\mathbf{H}}_{P2}^{H}$ represent the channel estimate between all BSs and the $k$th MS. We first compare the performance of the robust centralized and robust distributed algorithms in terms of the total sum power of all BSs. Fig. 5 shows that the robust centralized and distributed algorithms utilize the same total power. Then, we compare the performance of the proposed robust design and the non-robust design of [20][7] in terms of total sum power of all BSs which is also plotted

---
[6] Note that the feasible initial $\{w_k\}_{k=1}^{K}$ for problem $\mathcal{P}2$ can be obtained from the solution of $\mathcal{P}1$.

[7] For $\mathcal{P}2$, we have also shown in the latter paper that [17] and [20] have the same performance.

$$\widehat{\mathbf{H}}_{P2}^{H} = \begin{bmatrix} 0.2328 - 0.0868i & 0.1344 + 0.3848i & 0.2407 - 0.3118i & -0.2276 - 1.3829i \\ 1.6717 + 0.4976i & 0.5254 - 0.9034i & 0.0206 + 0.1318i & 1.1292 - 0.8314i \\ -0.0426 + 1.7262i & -0.6380 - 0.5663i & -0.2765 - 0.4716i & -0.2760 - 0.1349i \\ 0.3266 + 0.0882i & -0.3655 + 0.8997i & -0.5944 - 1.1556i & -0.2692 - 0.6244i \end{bmatrix} \quad (39)$$

in Fig. 5. This figure shows that the robust design utilizes more power than that of the non-robust design for all noise variances. Now, for the total power given in the latter figure, the AMSE of each user for both designs are plotted in Fig. 6. This figure shows that for all noise variances, the non-robust design does not satisfy the AMSE requirement but the proposed robust design ensures the AMSE requirement efficiently. To ensure the latter AMSE, however, the robust design utilizes more power than that of the non-robust design. This shows that Fig. 6 does not reveal the actual performance of the proposed robust design. Thus, for a fair comparison, we tune $\{\sigma_{ek}^2\}_{k=1}^K$ such that the total power utilized by the robust and non-robust designs are the same (or very close to each other)[8] and then we compare the performance of these two designs by their achieved AMSEs of each user. For the power requirement of Fig. 7, the AMSE of each user is plotted in Fig. 8. The latter figure shows that for the same total power, the robust design still outperforms the non-robust design.

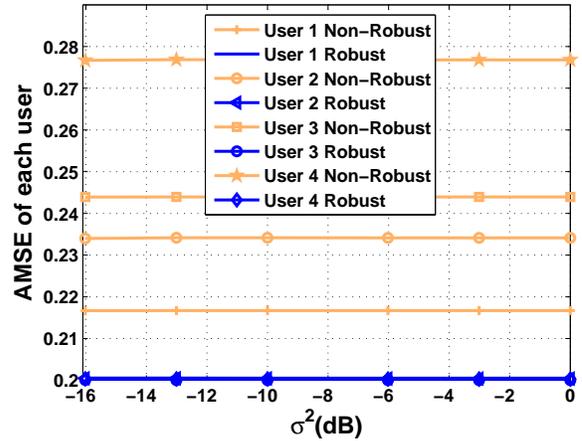

Fig. 6. Comparison of the AMSE achieved in each user for the robust design and non-robust design of [20].

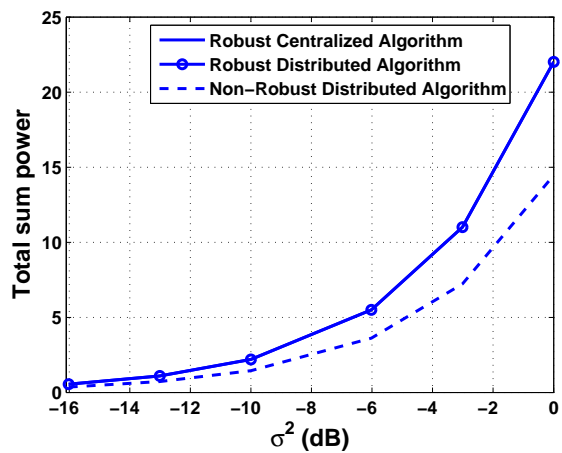

Fig. 5. Comparison of the robust centralized algorithm, distributed algorithm and the non-robust distributed algorithm of [20].

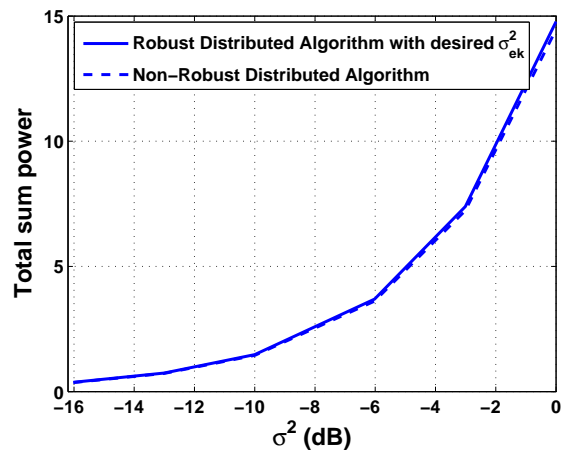

Fig. 7. The power utilized by the robust design and non-robust design of [20] with the desired error variances $\{\tilde{\sigma}_{ek}^2\}_{k=1}^K$.

For $\mathcal{P}2$, we have noticed that for each channel estimate different error variance (after numerical tuning) is required to get the same total sum powers for the robust and non-robust designs. However, the performance behavior exhibited in all channel estimates fits to that of $\widehat{\mathbf{H}}_{P2}^{H}$. The detail simulation results of this problem for the other channel estimates are omitted to reduce redundancy. Moreover, to see the effect of antenna correlation factor for $\mathcal{P}2$, we use the previously mentioned AMSE targets (i.e., $\{\varepsilon_k = 0.2\}_{k=1}^K$ ). With these

[8]In the robust design, for each noise level, we perform numerical search to get the appropriate $\{\sigma_{ek}^2\}_{k=1}^K$ that yields the same (or very close) total transmit power as that of the non-robust design. This task is termed as numerical tuning.

AMSE targets, when we increase the antenna correlation factor, we observe that the total power requirement of the whole network also increases. The simulation results which show this fact has not been included for conciseness.

### C. Convergence characteristics of **Algorithm I**

As we have mentioned in Section IV-A.1, the centralized algorithm to solve (7) has limited practical interest when the number of BSs and/or MSs are large. Moreover, in Section IV-A.2, the computational complexity of our distributed algorithm to solve (12) (once the complexity of (12) with **Algorithm I** is studied, the complexity of (7) with this algorithm is immediate) is provided for a single iteration (i.e., per iteration of **Algorithm I**). Thus, to show the computational advantage of our distributed algorithm compared to that of the centralized algorithm, the number of required iterations for convergence of **Algorithm I** ($I_1$) needs to be accessed for large scale networks. However, we are not able to compute $I_1$ analytically. Due to this, we examine the convergence characteristics of **Algorithm I** for a 19-cell hexagonal structure coordinated BS system as in [33]. Each BS is located at the center of its cell, whereas each MS is located randomly inside these 19-cells with uniform distribution. The propagation model between each BS and MS contains two components. One is the path loss component decaying with distance, and the other one is the Rayleigh fading random component which has a zero mean and unit variance. For this simulation, we use $\{\eta_k = 1, \sigma_{elk}^2 = 0.02, \rho_{blk} = 0.25, \forall l\}_{k=1}^K$ and all the other parameter settings are summarized as shown in Table I. For the channel realizations of these parameters, we examine the convergence characterstics of **Algorithm I** at different iterative stages of **Algorithm II** (i.e., with different $\{w_k\}_{k=1}^K$) as shown in Fig. 9. As can be seen from this figure, **Algorithm I** converges to an optimal solution in less than 10 iterations.

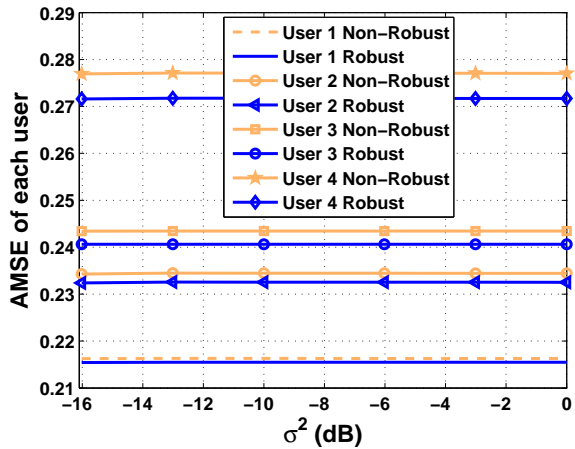

Fig. 8. Comparison of the AMSE achieved in each user for the robust design and non-robust design of [20] after tuning the error variance.

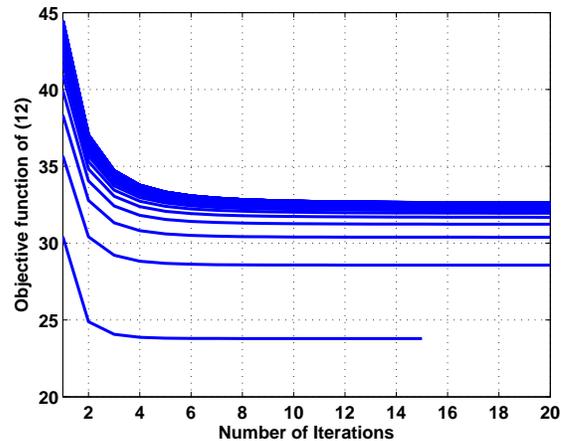

Fig. 9. Convergence characteristics of **Algorithm I** at different iterative stages of **Algorithm II**.

TABLE I
SIMULATION PARAMETERS FOR CONVERGENCE OF **Algorithm I**

| Number of BSs | 19 |
|---|---|
| Number of antennas at each BS | 2 |
| Transmit power of each BS antenna | 5W |
| Radius of each cell | 1.6km |
| Reference distance ($d_0$) | 1.6km |
| Path loss exponent | 3.8 |
| Mean path loss at $d_0$ | 134dB |
| Channel bandwidth | 5MHz |
| Receiver noise figure | 5dB |
| Receiver vertical antenna gain | 10.3dBi |
| Receiver temperature | 300K |
| SNR | 18dB |

### D. Overall computational complexity to solve (7)

In Section IV-A.1, we have presented the worst-case computational cost of IP methods to solve (7) centrally. However, in most practical problems, IP methods require less computational cost than that of their worst-case complexities. To the best of our knowledge, computing the exact computational complexity of IP methods for this problem requires immense effort and time. Hence, we believe that such a task is beyond the scope of our current work. However, we have carried out extensive simulations to compare the computational time of our proposed distributed algorithm with that of the centralized algorithm which uses IP method. In the following, we describe the simulation platform and methodology we have used, and discuss the results.

According to [34], MOSEC is a computationally efficient optimization package which uses IP methods to solve large-scale optimization problems. Moreover, for SOCP problems, MOSEC requires less computational time than that of SeDuMi, LOQO, SDPT3 and CPLEX [35], [36]. This motivates us to compare the computational time of **Algorithm I** with that of MOSEC to solve (7). Our Matlab codes were run on a personal computer with 1.6 GHz, 2GB dual core processor under Windows XP. For comparison between these two algorithms, we have used a coordinated BS system with $L = N/2$, $N = K$, $\{\rho_{blk} = 0.25, \sigma_{elk}^2 = 0.02, \sigma_k^2 = 0.1,\}_{k=1}^K$ and $\{p_n = 2\}_{n=1}^N$. It is assumed that problem (7) has been solved by a central controller with $K$ processors and all other parameters are taken as mentioned in the first paragraph of Section VI. Table II shows the amount of time required to solve (7) by **Algorithm I**[9] and MOSEC at different iterative stages of **Algorithm II** (i.e., for different $\{w_k\}_{k=1}^K$). As can be seen from Table II, our proposed distributed algorithm requires less computational time than that of MOSEC. From this table we can notice that our distributed algorithm has practical interest especially when

---

[9]To get the computational time of **Algorithm I** per processor, first we get the computational time of **Algorithm I** by assuming one processor (i.e., personal computer), then, we divide the latter computational time by K.

$K(N)$ is large.

TABLE II
COMPUTATIONAL TIME OF **Algorithm I** AND MOSEC FOR (7) (IN SECONDS)

| K | 4 | 10 | 20 | 30 | 40 |
|---|---|---|---|---|---|
| A: MOSEC | 0.0031 | 0.0213 | 0.2182 | 0.8833 | 2.3598 |
| B: $\frac{\textbf{Algorithm I}}{processor}$ | 0.0012 | 0.0020 | 0.0041 | 0.0086 | 0.0118 |
| A/B | 2.52 | 10.455 | 52.643 | 133.874 | 200.336 |

The convergence characteristics of **Algorithm III** and the overall computational complexity of (20) can be studied like in Sections VI-C and VI-D, respectively.

## VII. CONCLUSIONS

This paper considers the joint transceiver design for multiuser MISO systems with coordinated BSs where imperfect CSI is available at the BSs and MSs. By incorporating antenna correlation at the BSs and taking channel estimation errors into account, we solve two robust design problems. The problems are solved as follows. First, for fixed receivers, we propose centralized and novel computationally efficient distributed algorithms to jointly optimize the precoders of all users. The centralized algorithms employ the SOCP approach, whereas the distributed algorithms use the Lagrangian dual decomposition, modified matrix fractional minimization and an iterative method. Second, for fixed BS precoders, the receivers are updated by the MAMSE criterion. These two steps are repeated until convergence is achieved. Computer simulations demonstrate that our proposed distributed algorithms achieve the same performance as that of the centralized algorithms. Simulation results also verify the superior performance of the stochastic robust designs compared to that of the non-robust/naive designs.

## REFERENCES

[1] G. Caire and S. Shamai, "On the achievable throughput of a multi-antenna Gaussian broadcast channel," *IEEE Tran. Info. Theo*, vol. 49, no. 7, pp. 1691 – 1706, Jul. 2003.
[2] S. Serbetli and A. Yener, "Transceiver optimization for multiuser MIMO systems," *IEEE Tran. Sig. Proc.*, vol. 52, no. 1, pp. 214 – 226, Jan. 2004.
[3] D. P. Palomar and S. Verdu, "Gradient of mutual information in linear vector Gaussian channels," *IEEE Tran. Info. Theo.*, vol. 52, no. 1, pp. 141 – 154, Jan. 2006.
[4] D. P. Palomar, *A unified framework for communications through MIMO channels*, Ph.D. thesis, Technical University of Catalonia (UPC), Barcelona, Spain, 2003.
[5] S. Shi, M. Schubert, and H. Boche, "Downlink MMSE transceiver optimization for multiuser MIMO systems: Duality and sum-MSE minimization," *IEEE Tran. Sig. Proc.*, vol. 55, no. 11, pp. 5436 – 5446, Nov. 2007.
[6] S. Shi, M. Schubert, and H. Boche, "Rate optimization for multiuser MIMO systems with linear processing," *IEEE Tran. Sig. Proc.*, vol. 56, no. 8, pp. 4020 – 4030, Aug. 2008.
[7] R. Hunger, M. Joham, and W. Utschick, "On the MSE-duality of the broadcast channel and the multiple access channel," *IEEE Tran. Sig. Proc.*, vol. 57, no. 2, pp. 698 – 713, Feb. 2009.
[8] S. Shi, M. Schubert, and H. Boche, "Downlink MMSE transceiver optimization for multiuser MIMO systems: MMSE balancing," *IEEE Tran. Sig. Proc.*, vol. 56, no. 8, pp. 3702 – 3712, Aug. 2008.
[9] M. B. Shenouda and T. N. Davidson, "On the design of linear transceivers for multiuser systems with channel uncertainty," *IEEE Jour. Sel. Commun.*, vol. 26, no. 6, pp. 1015 – 1024, Aug. 2008.
[10] P. Ubaidulla and A. Chockalingam, "Robust joint precoder/receive filter designs for multiuser MIMO downlink," in *10th IEEE Workshop on Signal Processing Advances in Wireless Communications (SPAWC)*, Perugia, Italy, 21 – 24 Jun. 2009, pp. 136 – 140.
[11] T. Endeshaw, B. K. Chalise, and L. Vandendorpe, "MSE uplink-downlink duality of MIMO systems under imperfect CSI," in *3rd IEEE International Workshop on Computational Advances in Multi-Sensor Adaptive Processing (CAMSAP)*, Aruba, 13 – 16 Dec. 2009, pp. 384 – 387.
[12] T. E. Bogale, B. K. Chalise, and L. Vandendorpe, "Robust transceiver optimization for downlink multiuser MIMO systems," *IEEE Tran. Sig. Proc.*, vol. 59, no. 1, pp. 446 – 453, Jan. 2011.
[13] K. M. Karakayali, G. J. Foschini, and R. A. Valenzuela, "Network coordination for spectrally efficient communications in cellular systems," *IEEE. Tran. Wirel. Comm.*, vol. 13, no. 4, pp. 56 – 61, Aug. 2006.
[14] H. Dahrouj and W. Yu, "Coordinated beamforming for the multi-cell multi-antenna wireless system," *IEEE Tran. Wirel. Comm.*, vol. 9, no. 5, pp. 1748 – 1759, May 2010.
[15] E. Bjornson, R. Zakhour, D. Gesbert, and B. Ottersten, "Cooperative multicell precoding: Rate region characterization and distributed strategies with instantaneous and statistical CSI," *IEEE Tran. Sig. Proc.*, vol. 58, no. 8, pp. 4298 – 4310, Aug. 2010.
[16] E. Bjornson, R. Zakhour, D. Gesbert, and B. Ottersten, "Distributed multicell and multiantenna precoding: Characterization and performance evaluation," in *Proc. IEEE Global Telecommunications Conference (GLOBECOM)*, Honolulu, HI, USA, 30 Nov. – 4 Dec. 2009, pp. 1 – 6.
[17] S. Shi, M. Schubert, N. Vucic, and H. Boche, "MMSE optimization with per-base-station power constraints for network MIMO systems," in *Proc. IEEE International Conference on Communications (ICC)*, Beijing, China, 19 – 23 May 2008, pp. 4106 – 4110.
[18] S. Boyd and L. Vandenberghe, *Convex optimization*, Cambridge University Press, Cambridge, 2004.
[19] T. Tamaki, K. Seong, and J. M. Cioffi, "Downlink MIMO systems using cooperation among base stations in a slow fading channel," in *Proc. IEEE International Conference on Communications (ICC)*, Glasgow, UK, 24 – 28 Jun. 2007, pp. 4728 – 4733.
[20] T. E. Bogale, L. Vandendorpe, and B. K. Chalise, "MMSE transceiver design for coordinated base station systems: Distributive algorithm," in *44th Annual Asilomar Conference on Signals, Systems, and Computers*, Pacific Grove, CA, USA, 7 – 10 Nov. 2010.
[21] J. Luo, J. R. Zeidler, and S. Mclaughlin, "Performance analysis of compact antenna arrays with MRC in correlated Nakagami fading channels," *IEEE Trans. Veh. Technol*, vol. 50, pp. 267–277, 2001.
[22] B. K. Chalise, L. Haering, and A. Czylwik, "System level performance of UMTS-FDD with covariance transformation based DL beamforming," in *Proc. IEEE Global Telecommunications Conference (GLOBECOM)*, San Fransisco, USA, 1 – 5 Dec. 2003, pp. 133 – 137.
[23] D. Aszetly, *On antenna arrays in mobile communication systems: Fast fading and GSM base station receiver algorithm*, Ph.D. thesis, Royal Institute Technology, Stockholm, Sweden, 1996.
[24] T. Yoo, E. Yoon, and A. Goldsmith, "MIMO capacity with channel uncertainty: Does feedback help?," in *Proc. IEEE Global Telecommunications Conference (GLOBECOM)*, Dallas, USA, 29 Nov. – 3 Dec. 2004, vol. 1, pp. 96 – 100.
[25] W. Yu and T. Lan, "Transmitter optimization for the multi-antenna downlink with per-antenna power constraints," *IEEE Trans. Sig. Proc.*, vol. 55, no. 6, pp. 2646 – 2660, Jun. 2007.
[26] M. Ding, *Multiple-input multiple-output wireless system designs with imperfect channel knowledge*, Ph.D. thesis, Queens University Kingston, Ontario, Canada, 2008.
[27] D. Ding and S. D. Blostein, "MIMO minimum total MSE transceiver design with imperfect CSI at both ends," *IEEE Tran. Sig. Proc.*, vol. 57, no. 3, pp. 1141 – 1150, Mar. 2009.
[28] K. B. Petersen and M. S. Pedersen, "The matrix cookbook," Feb. 2008.
[29] M. S. Lobo, L. Vandenberghe, S. Boyd, and H. Lebret, "Applications of second-order cone programming," *Linear algebra and its applications*, vol. 284, pp. 193 – 228, 1998.
[30] D. Coppersmith and S. Winograd, "Matrix multiplication via arithmetic progressions," *Journal of Symbolic Computation*, vol. 9, pp. 251 – 280, 1990.
[31] None, "Cubic function: http://en.wikipedia.org/wiki/cubic function/," .
[32] S. L. Loyka, "Channel capacity of mimo architecture using the exponential correlation matrix," *IEEE Tran. Comm.*, vol. 5, no. 9, pp. 369 – 371, 2001.
[33] G. J. Foschini, K. Karakayali, and R. A. Valenzuela, "Coordinating multiple antenna cellular networks to achieve enormous spectral efficiency," *Communications, IEE Proceedings*, vol. 153, no. 4, 2006.

[34] MOSEK ApS, *The MOSEK optimization toolbox for MATLAB manual, Version 6.0 (Revision 103), http://www.mosek.com*, 2011.
[35] H. D. Mittelmann, "An independent benchmarking of SDP and SOCP solvers," *Mathimatical Programming*, vol. 95, no. 2, pp. 407 – 430, 2003.
[36] H. D. Mittelmann, *SOCP (second-order cone programming) Benchmark, http://plato.asu.edu/ftp/socp.html*, Feb. 2011.